\documentclass[moor]{informs1}              

\usepackage{natbib}
 \NatBibNumeric
 \def\bibfont{\small}%
 \bibpunct[, ]{[}{]}{,}{n}{}{,}%
\definecolor{DarkOrchid} {cmyk}{0.40,0.80,0.20,0}
\definecolor{VioletRed} {cmyk}{0,0.81,0,0}
\definecolor{Orange} {cmyk}{0,0.61,0.87,0}

\def\gph{\mathop{\rm gph\,}}

\def\cone{\mathop{\rm cone}}
\def\conv{\mathop{\rm conv}}
\def\cl{\mathop{\rm cl}}

\def\int{\mathop{\rm Int }}
\def\im{\mathop{\rm Im}}
\def\dom{\mathop{\rm Dom\,}}
\def\reg {\mathop{\rm reg\,}}

\def\N{{\mathbb{N}}}

\def\R{{\mathbb{R}}}

\def\TheoremsNumberedThrough{%
\theoremstyle{TH}%
\newtheorem{theorem}{Theorem}
\newtheorem{lemma}{Lemma}
\newtheorem{proposition}{Proposition}
\newtheorem{corollary}{Corollary}

\theoremstyle{EX}
\newtheorem{remark}{Remark}

\newtheorem{definition}{Definition}

}
\usepackage[colorlinks=true,breaklinks=true,bookmarks=true,urlcolor=blue,
     citecolor=blue,linkcolor=blue,bookmarksopen=false,draft=false]{hyperref}

\def\EMAIL#1{\href{mailto:#1}{#1}}

\TheoremsNumberedThrough     

\EquationsNumberedThrough    

\begin{document}


 \RUNAUTHOR{Huynh and Th\'era}

\RUNTITLE{Directional metric regularity of multifunctions}

 \TITLE{Directional metric regularity of multifunctions}

\ARTICLEAUTHORS{%
\AUTHOR{Huynh Van Ngai }
\AFF{Department of Mathematics, University of
Quy Nhon, 170 An Duong Vuong, Quy Nhon, Viet Nam \\\EMAIL{nghiakhiem@yahoo.com}
}
\AUTHOR{Michel Th\'era}
\AFF{Laboratoire XLIM, Universit\'e de Limoges\\\EMAIL{michel.thera@unilim.fr}
}
} 

\ABSTRACT{%
In this paper, we study relative metric regularity of set-valued mappings with emphasis on  directional metric regularity. We establish characterizations of  relative metric regularity without assuming  the completeness of  the image spaces, by using the relative lower semicontinuous envelopes of the distance functions to set-valued mappings. We then apply these characterizations to establish a coderivative type criterion for  directional metric regularity as well as for the robustness of  metric regularity.}

\MSCCLASS{49J52, 49J53, 90C30.}
\ORMSCLASS{Primary: Variational Analysis ; secondary:  Optimisation}

\maketitle

%
\section{Introduction}\label{intro}

Set-valued mappings represent the most developed class of objects studied in the framework of variational analysis. Various types of set-valued mappings arise in a considerable number of models ranging from mathematical programs, through game theory and to control and design problems.
The most well known and widely used regularity property of set-valued mappings is that of \emph{metric regularity} \cite{BorweinZhuBook, DmiKru09.1,
Io00, KlaKum02, Pen89, RockBook2002, Borisbook1, DonRoc09, DmiKru08.1,NgaiThe08,
Gia05,Kum99}.
The term ``metric regularity'' was coined by Borwein  \& Zhuang \cite{BorZhuang88}.
Metric regularity  or its equivalent notions  (openness or covering at a linear rate or Aubin property of the inverse)   is  a central concept in modern variational analysis. In particular, this property is used as a key  ingredient in  investigating the behavior of the solution set of generalized equations associated to set-valued mappings.

According to the long history of metric regularity there is an abundant  literature on conditions ensuring this property. The roots of this notion go back to the classical \textit{Banach Open Mapping Theorem}  (see, for instance \cite{conway},  Theorem III.12.1) and its subsequent generalization to nonlinear mappings known as \textit{Lyusternik and Graves Theorem}(\cite{Lyusternik, RefGr}, see also  \cite{RefDMO80, dontchev}).
  For a detailed account  on results on   metric regularity as well as on  its various applications,  we refer the reader    to   basic monographs  and  references, \cite{Refaze, Aze06, RefBonS, BD, BorweinZhuBook, BorZhu96,BorZhuang88, Clarke83, RefCom, DmiKru08.1, DonLewRock03, RefIo3, Io00, RefIo2, IofOut08, RefJT, JT1, Lyusternik, Borisbook1, Borisbook2, RefMorS, RefMorW, RefNT1, RefNT3, Outrata98, Pen89, penot2013, robinson80}, 
as well as to  the references given therein. 
\vskip 0.2cm
 Let $X$ and $Y$ be metric spaces endowed with  metrics both denoted by $d(\cdot,\cdot).$ The open ball with center $x$  and radius $r>0$ is denoted by $B(x,r).   $ For a given subset $\Omega$ of $X$,  we denote by $\conv \Omega$ and  $\cone \Omega$ the convex hull of $\Omega$ and the conical
convex hull of $\Omega$, respectively.  We also  use the symbols $w$ and $w^\star$ to indicate the weak and the weak$^\star$
topology, respectively, and w- lim and w$^\star$- lim represent the weak and the weak$^\star$ topological limits, respectively. $\int \Omega$  and
$\cl  \Omega$ are the interior and the closure of $\Omega$  with respect to the norm topology, respectively; $\cl ^w\Omega $ stands for the closure in
the weak topology  and $\cl ^{w^\star }\Omega$  is the closure in the weak$^\star$ topology of a given subset $\Omega \subset X^\star$ in the dual space. We also make
use of the property,  that for convex sets, the norm and the weak closures coincide. Finally given a mapping $f:X\to Y$ we note $\im f$  for  the range   $f.$
 
Recall that   a set-valued (multivalued) mappping $F: X\rightrightarrows Y$  is a mapping which assigns to every $x\in X$ a subset (possibly empty) $F(x)$ of   $Y$. As usual, we use the notation  $\gph F:= \{(x,y)\in X\times Y\ :\ y\in F(x)\}$  for  the graph of $F$,
$\dom F:=\{x\in X : F(x)\neq \emptyset \}$ for the domain  of $F$
and $F^{-1}:Y\rightrightarrows X$  for the  inverse  of $F$. This
inverse (which always exists) is defined by $F^{-1}(y):=\{x\in X\ :\
y\in F(x)\}, \,\, y\in Y$ and satisfies  $$(x,y)\in \gph F\;\iff\;
(y,x)\in \gph F^{-1} .$$ 
If  $C$ is a subset of $X$, we    use  the standard notation  $d(x, C) = \inf_{z \in C} d(x,z)$,  with the convention that   $d(x, S) = +\infty$ whenever $C$ is empty. We recall  that a   multifunction $F$ is  {\it metrically regular} at  $( x_0,y_0)\in\gph F$   with  modulus  $\tau>0$ if there
exists a neighborhood $B( (x_0,y_0),\varepsilon)$ of $( x_0, y_0)$  such that
\begin{equation}\label{Regular}
d(x,F^{-1}(y))\leq \tau d(y,F(x))\quad\mbox{for all} \quad (x,y)\in
B( (x_0,y_0),\varepsilon).
\end{equation}
The infimum  of all moduli $\tau$   satisfying relation (\ref{Regular})  is denoted by  $\reg F(\bar x,\bar y)$\; (\cite{DonLewRock03}).
 In the case for example  of a set-valued mapping $F$ with  a closed and convex graph,   the Robinson-Ursescu Theorem   (\cite{robinson73} and  \cite{urcescu}),    says  that  $F$ is metrically regular at $(x_0,y_0)$,  if and only if   $y_0$ is an interior point to the range of $F$, i.e., to $\dom F^{-1}$.
 \vskip 0.2cm
Recently, several generalized or weaker versions of metric regularity (restricted metric regularity \cite{{RefMorW}}, calmness, subregularity) have been considered. Especially, Ioffe (\cite{RefIoffe})  introduced and studied a natural extension of  metric regularity called  ``\textit{relative metric regularity}" which covers almost  every notions of metric regularity  given in the literature. Roughly speaking, a mapping $F$ is relatively metrically  regular relative to some subset $V\subseteq X\times Y$  if the metric regularity  property is satisfied at points belonging to $V$ and near the reference point. An important special case of this relative metric regularity concept is the notion of  directional metric regularity introduced and studied by Arutyunov, Avakov and Izmailov in \cite{ArAvIz07}.
This directional metric regularity is an extension of an earlier concept used by Bonnans \& Shapiro (\cite{RefBonS}) to study sensitivity analysis. 
\vskip 0.2cm
Our  main  objective in this paper is to use the theory of  error bounds  to study
directional metric regularity of multifunctions.  We develop  the method used  by the authors in  \cite{RefNT3, NgaiThe08, NTT} to characterize   relative metric regularity by using  global/local slopes of a suitable lower semicontinuous envelope type  of the distance function to the images of set-valued mappings. A particular advantage of this approach is to avoid the  completeness of the image space that is not really necessary in some important situations. These established characterizations permit   to derive   coderivative conditions as well as   stability results for  directional metric regularity.  
\vskip 0.2cm
The remainder of this paper is organized as follows.
In Section 2, we prove  characterizations  
of   relative metric regularity of closed multifunctions on metric spaces by using global/local strong slopes of suitable relative  semicontinuous envelope of distance functions to the images of set-valued mappings.  Based on 
these characterizations, we derive in Section 3 coderivative criteria  ensuring
directional metric regularity. In the final section, results on the perturbation stability of  directional metric regularity are reported.
 \section{Characterizations of relative metric regularity}
 Let  $X$  be a metric space. Let  $f : X \rightarrow \mathbb{R} \cup
\{+\infty\}$  be a given function.  As usual, $\mbox{dom}f := \{x
\in X : f(x) < +\infty\}$ denotes the domain of $f$. We set
\begin{equation}\label{1}
S:=\{x\in X:\quad f(x)\leq 0\}.
\end{equation}
We use the symbol $[f(x)]_+$ to denote $\max(f(x),0).$ We shall say
that the system (\ref{1}) admits  {\it a global  error bound}
 if there exists  a real $c > 0$  such that \begin{equation}\label{2}
d(x, S) \leq c\big[f(x)]_+ \quad\mbox{ for all}\quad  x\in X.
\end{equation}
For $x_0 \in S$, we shall say that the system  (\ref{1}) has a \textit{ local 
error bound } at $x_0$, when there exist reals  $c>0$ and $
\varepsilon > 0$  such that  relation (\ref{2}) is satisfied for all
$x$ around  $x_0$, i.e.,  in an open  ball  $B(x_0,\varepsilon)$
with center $x_0$ and radius $\varepsilon.$
\vskip 0.2cm 
Since the first error bound result due to  Hoffman (\cite{Hof52}), numerous characterizations and criteria for error bounds in terms of various derivative-like objects have been established.
\cite{FabHenKruOut10,FabHenKruOut12}. Stability and some other properties of error bounds are examined in \cite{NgaiKruThe10,NgaiThe05,NgaiThe08,NgaiThe09,KruNgaiThe10}.  Several
conditions using subdifferential operators or directional
derivatives and ensuring the error bound property  in Banach spaces have been
established, for example,  in   \cite{BJH, RefJ, RefNT3}. Error bounds  have been also used in sensitivity analysis of linear programming/linear complementarity problem and also as termination criteria for descent algorithms.
Recently, Az\'e
\cite{Refaze}, Az\'e \& Corvellec \cite{RefAC2} have used the
so-called strong slope introduced by De Giorgi, Marino \& Tosques
in \cite{RefDMT} to prove criteria  for error bounds in complete
metric spaces. 
\vskip 0.2cm
Recall from \cite{Io00}\footnote {A. Ioffe was at your knowledge the first one to  start using slopes and to advertise them in the optimization community.}  that the
local and  global strong slopes $\vert\nabla f\vert (x);$ $|\Gamma f|(x)$ of a function $f$ at $x\in\mbox{dom}f$ are the quantities defined by
\begin{equation}\label{Slopes}
\vert\nabla f\vert (x)=\limsup_{y\to x, \;y\neq x}\frac{[f(x)-f(y)]_+}{d(x,y)};\quad\vert\Gamma f\vert (x)=\sup_{y\not= x}\frac{[f(x)-f(y)]_+}{d(x,y)}
\end{equation}
For $x\notin \mbox{dom}f,$  we set  $\vert\nabla f\vert
(x)=|\Gamma f|(x)=+\infty.$   When $f$ takes only negative values it coincides with $$|\nabla{f}|^{\diamond}(x):=\sup_{y\neq{x}}
\frac{[f(x)-f_+(y)]_+}{d(y,{x})}$$
as defined in \cite{NgaKruThe12}. 
As pointed out by Ioffe \cite[Poposition 3.8]{ioffe2013}, when $f$ is a convex function  defined on a Banach space, then  
$$\vert\nabla f\vert (x)= \sup_{\Vert h\Vert\leq 1} (-\vert f^\prime(x;h))= d(0, \partial f(x).$$
Trivially, one has $\vert\nabla f\vert (x)\le |\Gamma f|(x),$ for all $x\in X.$
\vskip 0.2cm In the sequel, we will need the
following result established by Ngai \& Th\'era (\cite {NgaiThe08}), which gives an
estimation via the global strong slope for the distance $d(\bar{x},S)$ from a given point
$\bar{x}$ outside of $S$ to the set $S$ in complete metric spaces.
\begin{theorem}\label{T1}
Let $X$ be a complete metric space  and let $f : X\to
\R\cup\{+\infty\}$ be a lower semicontinuous function and
$\bar{x}\notin S$.  Then, setting
\begin{equation}\label{m}
m(\bar{x}):=\inf\left\{|\Gamma f|(x):\quad  
d(x, \bar{x})   < d(\bar{ x},S),\;
f(x) \leq f( \bar {x})
\right\},\end{equation}
one has
\begin{equation}\label{Estim}
m(\bar{x})d(\bar{x},S)\leq f(\bar{x}).
\end{equation}
\end{theorem}
\vskip 0.5cm
Let $X$ be a  metric space and let $Y$ be a normed linear space. Consider a multifunction
$F:X\rightrightarrows Y$. Let us recall from (\cite{ArAvIz07}) the definition of directional metric regularity.
\begin{definition}\label{Def-Direct-Regular} Let $F:X\rightrightarrows Y$ be a multifunction. Let $(x_0,y_0)\in\gph F$ and $\bar y\in Y$ be given. $F$ is said to be directionally metrically regular at $(x_0,y_0)$ in the direction $\bar y$ with a modulus $\tau>0$ if there exist $\varepsilon>0,$ $\delta>0$
such that 
\begin{equation}\label{DMR} 
d(x, F^{-1}(y))\le\tau d(y, F(x)), 
\end{equation}
for all $(x,y)\in B(x_0,\varepsilon)\times B(y_0,\varepsilon)$ satisfying 
$$ d(y, F(x))<\varepsilon\quad\mbox{and}\quad y\in F(x)+\cone B(\bar y,\delta) .$$
Here, $\cone B(\bar y,\delta)$ stands for the conic hull of $B(\bar y,\delta)$, i.e., $  \cone B(\bar y,\delta)= \displaystyle\bigcup_{\lambda \geq 0}\lambda B(\bar y,\delta)$.
The infinum of all moduli  $\tau$  in relation (\ref{DMR}) is called the \textit{exact modulus} of the metric regularity at $(x_0,y_0)$ in direction $\bar y,$ and is denoted by $\reg_{\bar y} F(x_0,y_0).$
\end{definition}
\vskip 0.2cm
Note that if $F$ is metrically regular at $(x_0,y_0)\in\gph F$,  then $F$ is directionally metrically regular in all directions $\bar y\in Y.$  When $\|\bar y\|<\delta,$ then   directionally metric regularity coincides with the usual metric regularity. 
The notion of   directionally metric regularity is a special case of \textit{ metric regularity relative to a set } $V$ with $\gph F\subseteq V\subseteq X\times Y,$ introduced by Ioffe (\cite{RefIoffe}). For $y\in Y,$ $x\in X,$ denote by $V_y:=\left\{x\in X:\; (x,y)\in V\right\},$  $V_x:=\left\{y\in Y:\; (x,y)\in V\right\},$ and $\cl V_y$  for the closure of $V_y.$
\begin{definition}\label{Def-Rela-Regular} Let $X, Y$ be metric spaces. Let $F:X\rightrightarrows Y$ be a multifunction and  let $(x_0,y_0)\in\gph F$ and  fix  a subset $V\subseteq X\times Y$. $F$ is said to be metrically regular relative to $V$ at $(x_0,y_0)$ with a constant $\tau>0$ if there exist $\varepsilon>0$ 
such that 
\begin{equation}\label{MRre@} 
d(x, F^{-1}(y)\cap \cl V_y)\le\tau d(y, F(x)), \quad\mbox{for all}\;\; (x,y)\in B((x_0,y_0), \varepsilon)\cap V,\;d(y,F(x))<\varepsilon.
\end{equation}
The infinum of all moduli  $\tau$ is called the exact modulus of  metric regularity at $(x_0,y_0)$ relative to $V,$ and denoted by $\reg_{V} F(x_0,y_0).$

\end{definition}
\vskip 0.2cm
 In papers \cite{NgaiThe08, NTT},  the lower semicontinuous envelope
$x\mapsto \varphi(x,y)$  of the function $x\mapsto d(y,
F(x))$ for $y\in Y$ i.e., 
$$\varphi(x,y):=\liminf_{u\to x}d(y,F(u)),$$
has been used to characterize   metric regularity of $F.$ Along with the relative metric regularity, we define for each $y\in Y$  the lower semicontinous envelope of the functions $x\mapsto d(y, F(x))$ relative to a set $V\subseteq X\times Y$ by setting
\begin{equation}
\varphi_V(x,y):=\left\{\begin{array}{lll}&\liminf_{V_y\owns u\to x}d(y, F(u))&\quad\mbox{if}\quad x\in\cl V_y\\
&+\infty&\quad\mbox{otherwise.}
\end{array}\right.
\end{equation}
Obviously, for each $y\in Y,$ the function $\varphi_V(\cdot,y)$ is lower semicontinuous.
\vskip 0.2cm
For  a given $\bar y\in Y, $   directionally metric regularity  in a given direction $\bar y$ is exactly metric regularity relative to $V(\bar y,\delta)$ (for some $\delta>0$): 
\begin{equation}\label{V}
V(\bar y,\delta):=\{(x,y):\; y\in F(x)+\cone B(\bar y,\delta)\}.
\end{equation}
In this case, the lower semicontinuous envelope function relative to $V(\delta,\bar y)$ is denoted simply by
\begin{equation} \label{envelope}
\varphi_\delta(x,y):=\varphi_ {V(\delta,\bar y)}(x,y).
\end{equation}
The following proposition permits to transfer equivalently  relative metric regularity of $F$ to the error bound property of the function $\varphi_V.$ 
\begin{proposition} \label{trans}
Let $F:X\rightrightarrows Y$ be a closed multifunction  (i.e., its graph is closed)  and let $(x_0,y_0)\in\gph F.$ For $V\subseteq X\times Y,$ the following statements holds.

\item {(i)} For all $y\in Y,$ one has
$$ F^{-1}(y)\cap \cl V_y=\{x\in X:\; \varphi_V(x,y)=0\};$$

\item{(ii)} $F$ is metrically regular relative to $V$ at $(x_0,y_0)$ with a modulus $\tau>0$ if and only if there exists $\varepsilon>0$ such that
$$ d(x,F^{-1}(y)\cap \cl V_y)\le \tau\varphi_V(x,y)\quad\mbox{for all}\;\; (x,y)\in B(x_0,\varepsilon)\times B(y_0,\varepsilon)\;\;\mbox{with}\;\; d(y, F(x))<\varepsilon.$$
\end{proposition}
\vskip 0.2cm
\proof{Proof.}  The proof follows  straightforwardly from the definition.\hfill\Halmos
\vskip 0.2cm
We establish in the next  theorem   characterizations of  relative metric regularity by using the local/global  strong slopes. 
\begin{theorem}\label{Char-strong-slope}

Let $X$ be a complete metric space and $Y$ be a metric space. Let
  $F:X\rightrightarrows Y$  is a closed multifunction  and let $(x_0,y_0)\in \gph F$ and let $ V\subseteq X\times Y.$   For  a given $\tau\in (0,+\infty),$ consider  the following statements.

\item{(i)}  $F$ is metrically regular relative to $V$ at $(x_0,y_0);$

\item{(ii)} There exists $\delta>0$ such that
\begin{equation}\label{global-char$}
|\Gamma\varphi_V(\cdot,y)|(x)\ge\tau^{-1}\;\;\mbox{for all}\;  (x,y)\in \big(B(x_0,\delta)\times B(y_0,\delta)\big),\;x\in\cl V_y \;\;\mbox{with}\;d(y,F(x))\in (0,\delta);
\end{equation}
\item{(iii)} There exists $\delta>0$ such that
\begin{equation}\label{local-char@@@}
|\nabla\varphi_V(\cdot,y)|(x)\ge\tau^{-1}\;\;\mbox{for all}\;  (x,y)\in \big(B(x_0,\delta)\times B(y_0,\delta)\big),\;x\in\cl V_y\;\;\mbox{with}\;d(y,F(x))\in (0,\delta);
\end{equation}

\item {(iv)} There exist $\delta>0$ such that
\begin{equation}\label{local-char@@@@}
|\nabla\varphi_V(\cdot,y)|(x)\ge\tau^{-1}\;\;\mbox{for all}\;  (x,y)\in \big(B(x_0,\delta)\times B(y_0,\delta)\big)\cap V\;\;\mbox{with}\;d(y,F(x))\in (0,\delta).
\end{equation}
Then, $(i)\Leftrightarrow(ii)\Leftarrow (iii)\Rightarrow (iv).$  In addition, if $Y$ is a normed linear space; $\gph F\subseteq V;$ $V_x$ is convex for any $x$ near $x_0$ and $V_y$ is open for $y$ near $y_0,$ then $(i)\Rightarrow (iv)$.

\end{theorem}
  \proof{Proof.}  The implications $(iii)\Rightarrow (ii)$ and $(iii)\Rightarrow (iv)$  are obvious. For $(i)\Rightarrow (iii),$ assume that $F$ is metrically regular relative  to $V$ at $(x_0,y_0)$ with modulus $\tau>0$.Then, there is $\delta>0$ such that
 $$d(x,F^{-1}(y))\leq \tau \varphi_V(x,y)\;\;\forall (x,y)\in B(x_0,\delta)\times B(y_0,\delta).$$
Let $(x,y)\in B(x_0,\delta)\times B(y_0,\delta)$ with  $\varphi_V(x,y)\in (0,+\infty)$  be given.  For any $\varepsilon>0,$ we can find $u\in F^{-1}(y)$ satisfying
$$ d(x,u)\le (\tau+\varepsilon)\varphi_V(x,y).$$
Then, $\varphi_V(u,y)=0,$
 $u\not=x$ and therefore,
$$ |\Gamma\varphi_V(\cdot,y)|(x)\ge \frac{\varphi_V(x,y)-\varphi_V(u,y)}{d(x,u)} =\frac{\varphi_V(x,y)}{d(x,u)}\ge (\tau+\varepsilon)^{-1}.$$
Since $\varepsilon>0$ is arbitrary, $(ii)$ holds.

Let us prove $(ii)\Rightarrow(i).$  
Suppose that (\ref{global-char$}) is satisfied for $\delta>0$.
Let $\varepsilon\in (0,\tau/2)$ be given, and let
$$\alpha:=\min\Big\{\delta/2,\frac{\delta}{2(\tau+\varepsilon)},\delta\tau\Big\}.$$

Let $(x,y)\in B((x_0,y_0),\alpha)$ with $x\in\cl V_y,$ $d(y, F(x))<\alpha$ be given. Then,
$$\varphi_V(x,y)<\inf _{u\in X}\varphi_V(u,y)+\alpha.$$
By virtue of the Ekeland variational principle \cite{RefE}  applied to
the function $x\mapsto\varphi_V(x,y)$ on $X,$ we can find $z\in X$
satisfying $d(x,z)\leq \alpha(\tau+2\varepsilon)$ and
$\varphi_V(z,y)\leq\varphi_V(x,y) (<\alpha)$ such that
$$\varphi_V(z,y)\leq\varphi_V(u,y)+\frac{1}{\tau+2\varepsilon}d(u,z)\quad\mbox{for all}\; u\in X.$$
Consequently, $z\in B(x_0,\delta)$ and
$$\varphi_V(z,y)-\varphi_V(u,y)\leq \frac{d(z,u)}{\tau+2\varepsilon}\leq\frac{d(z,u)}{\tau+\varepsilon}\quad \mbox{for all}\; u\in X.$$
Therefore, by relation (\ref{global-char$}), we must have $z\in F^{-1}(y)\cap \cl V_y.$ Consequently,
\begin{equation}\label{Nonempty}
B(x_0,2\alpha\tau)\cap F^{-1}(y)\cap \cl V_y\not=\emptyset.
\end{equation}
 Then for any $z\in
X$ with $d(x,z)<d(x,F^{-1}(y));$ $\varphi_V(z,y)\leq \varphi_V(x,y),$ one has:
$$d(z,x_0)\leq d(z,x)+d(x,x_0)\leq d(x_0,F^{-1}(y))+2d(x,x_0)<2\alpha\tau+2\alpha\le\delta.$$
Thus, $z\in B(x_0,\delta)$ and $z\notin F^{-1}(y).$ Therefore, according to (\ref{global-char$}), one has
\begin{equation}
m(x):=\inf\left\{|\Gamma\varphi_(\cdot,y)|(z):\quad
\begin{array}{ll}
&d(z, x)   <d(x, F^{-1}(y))\\
&\varphi_V(z,y) \leq \varphi_V(x,y)
\end{array}\right\}>\frac{1}{\tau+\varepsilon}\notag.
\end{equation}
By virtue of Theorem \ref{T1} and as $\varepsilon>0$ is arbitrarily
small, we obtain
\begin{equation}\label{Inq}
d(x,F^{-1}(y)\cap V_y)\leq \tau \varphi_{V}(x,y),\notag
\end{equation}
which proves $(ii)\Rightarrow (i).$
\vskip 0.2cm
To conclude  the proof of the theorem, we need to show $(i)\Rightarrow (iv)$ provided $Y$ is a normed linear space; $V_y$ is open for $y$ near $y_0$ and $V_x$ is a convex set for $x$ near $x_0.$
 Let $\delta\in (0,1)$  be such that $V_y$ is open for all $y\in B(y_0,\delta);$ $V_x$ is convex for all $x\in B(x_0,\delta)$ and that
$$d(x,F^{-1}(y))\le \tau d(y,F(x))\quad\forall (x,y)\in \big(B(\bar x,2\delta)\times B(\bar y,2\delta)\big)\cap V.$$ Let $(x,y)\in (B(\bar x,\delta)\times
B(\bar y,\delta))\cap V$ be given with $F(x)\not=\emptyset;$ $y\notin F(x);$ $d (y,F(x))<\delta.$
For any $\varepsilon\in(0,\delta/2),$ pick
$u_\varepsilon\in B(x_0,\varepsilon)\cap V_y$ such that
\begin{equation}\label{Con-Seq} d(u_\varepsilon,x)<\varepsilon^2\varphi_V(x,y);\; d(y,F(u_\varepsilon))\leq (1+\varepsilon^2)^{1/2}\varphi_V(x,y).
\end{equation}
Take $y_{\varepsilon}\in F(u_\varepsilon)$ such that
$$d(y,F(u_\varepsilon))\leq\|y-y_\varepsilon\|<(1+\varepsilon^2)^{1/2}d(y,F(u_\varepsilon)).$$
Then, since $u\in B(x_0,\delta),$ by the convexity of $V_u,$ 
  $$ (u, z_\varepsilon) \in V\quad\mbox{with}\;\; z_\varepsilon:=\varepsilon y+(1-\varepsilon)y_\varepsilon.$$
Furthemore,
$$\|y-z_{\varepsilon}\|=(1-\varepsilon)\|y-y_{\varepsilon}\|<(1-\varepsilon)(1+\varepsilon^2)^{1/2}d(y,F(u_\varepsilon))<d(y, F(u_\varepsilon))<(1+\varepsilon)\varphi_V(x,y).$$
Therefore, $z_{\varepsilon}\notin F(u_\varepsilon)$ and $\|z_{\varepsilon}- y_0\|\leq
\|y-y_0\|+\|y-z_{\varepsilon}\|<2\delta.$ Hence,
we can select $x_{\varepsilon}\in F^{-1}(z_{\varepsilon})$ such that
\begin{equation}\begin{array}{ll}\label{Estim dist}
d(u_\varepsilon,x_{\varepsilon})<(1+\varepsilon)d(u_\varepsilon,F^{-1}(z_{\varepsilon}))&\leq
(1+\varepsilon)\tau d(z_{\varepsilon},F(u_\varepsilon))\\ &\leq
(1+\varepsilon)\varepsilon\tau\|y-y_{\varepsilon}\|.
\end{array}
\end{equation}
Consequently,  $\lim_{\varepsilon\to 0^+}d(x,x_{\varepsilon})=0.$ Hence, for $\varepsilon>0$ sufficiently small, $x_\varepsilon\in V_y,$ and 
 one has the following estimation 
\begin{equation}\label{Estim phi}
\begin{array}{ll}
\varphi_{V}(x,y)-\varphi_{V}(x_{\varepsilon},y)&\geq
\displaystyle\frac{1}{(1+\varepsilon^2)^{1/2}}d(y,F(u_\varepsilon))-d(y,F(x_\varepsilon))\\
&>
\displaystyle\left(\frac{1}{1+\varepsilon^2}-(1-\varepsilon)\right)\|y-y_{\varepsilon}\|\\&=
\displaystyle\frac{\varepsilon-\varepsilon^2+\varepsilon^3}{1+\varepsilon^2}\|y-y_{\varepsilon}\|.
\end{array}
\end{equation}
By combining  this relation and relations (\ref{Con-Seq}), (\ref{Estim dist}), one obtains
$$\frac{\varphi_{V}(x,y)-\varphi_{V}(x_{\varepsilon},y)}{d(x,x_{\varepsilon})}\geq \displaystyle\frac{\varphi_{V}(x,y)-\varphi_{V}(x_{\varepsilon},y)}{d(x,u_\varepsilon)+d(u_\varepsilon,x_{\varepsilon})}>
\frac{(\varepsilon-\varepsilon^2+\varepsilon^3)\|y-y_\varepsilon\|}{(1+\varepsilon^2)(\varepsilon^2\varphi_V(x,y)+(1+\varepsilon)\varepsilon\tau\|y-y_\varepsilon\|)}.$$
Since $\lim_{\varepsilon\to 0^+}\|y-y_\varepsilon\|=\lim_{\varepsilon\to0^+}d(y,F(u_\varepsilon))=\varphi_V(x,y)>0,$ then
$$|\nabla\varphi_V(\cdot,y)|(x)\geq \liminf_{\varepsilon\to0^+}\frac{\varphi_V(x,y)-\varphi_V(x_\varepsilon,y)}{d(x,x_{\varepsilon})}\geq \tau^{-1},$$
which completes the proof.\hfill$\Box$
\vskip 0.5cm  Theorem \ref{Char-strong-slope} yields the following exact formula
for the relative metric regularity.
\begin{corollary} Let $X$ be a complete metric space and let $Y$ be a metric space and let $V\subseteq X\times Y$ with $\gph F\subseteq V.$
Suppose that the multifunction $F:X\rightrightarrows Y$ is  closed
 and $(x_0,y_0)\in\mbox{gph}F$. Then, one has

$$1/\mbox{reg}_VF(x_0,y_0)=\liminf_{
\substack{(x,y) \overset {\varphi}{\rightarrow}(x_0,y_0)\\
y\notin F(x)\\ x\in\cl V_y}}|\Gamma\varphi_V(\cdot,y)|(x).$$
\vskip 0.2cm
Moreover, if  in addition, $Y$ is a normed linear space; $V_x$ is convex for any $x$ near $x_0$ and $V_y$ is open for $y$ near $y_0,$ then
$$1/\text{reg}_VF(x_0,y_0)\le\liminf_{
\substack{
(x,y)\overset {\varphi,V} {\rightarrow}(x_0,y_0)\\ y\notin F(x)}}|\nabla\varphi_V(\cdot,y)|(x).$$
The notation $(x,y)\overset {\varphi,V} {\rightarrow}(x_0, y_0)$  means that
$(x,y)\to (x_0, y_0)$ with $\varphi(x,y)\to 0$ and $(x,y)\in V.$
\end{corollary}
\vskip 0.2cm\proof{Proof.}  
It follows directly from Theorem \ref{Char-strong-slope}. \hfill$\Box$\vskip 0.5cm
We next introduce the partial notion of relative metric regularity for a parametric set-valued mapping.
Let $X,Y$ be metric spaces and let $P$ be a topological space. Given a set-valued mapping 
$F:X\times P\rightrightarrows Y$, we consider the implicit multifunction: $S: Y\times P\rightrightarrows Y$ defined by
\begin{equation}\label{Imp Multi}
S(y,p):=\{x\in X:\; \; y\in F(x,p)\}.
\end{equation}
Let $x_0\in S(y_0,p_0)$ and a set $V$: $\gph F\subseteq V\subseteq X\times P\times Y$ be given.
 \begin{definition}\label{Def-Uni-Direct-Regular} The set-valued mapping $F$ is said to be metrically regular uniformly in $p$ relatively  to $V$ at $((x_0,p_0), y_0)$ with a modulus $\tau>0$,  if there exist $\varepsilon>0$  and a neighborhood $W$ of $p_0$
such that 
\begin{equation}\label{MRre} 
d(x, S(y,p))\le\tau d(y, F(x,p)), \quad\mbox{for all}\;\; (x,y,p)\in \big(B((x_0,y_0)\varepsilon)\times W\big)\cap V;\;d(y, F(x,p))<\varepsilon.
\end{equation}
The infinum of all moduli  $\tau$ is called the exact modulus of the metric regularity of $F$  uniformly in $p$ at $ (x_0,y_0)$ relative to $V$  and is denoted by $\reg_{V} F(x_0,p_0,y_0).$

\end{definition}
\vskip 0.2cm 
Denote by
$$ V_{(x,p)}:=\{y\in Y:\;\; (x,p,y)\in V\},\;\;\mbox{for}\;(x,p)\in X\times P;$$  
$$ V_{(y,p)}:=\{x\in X:\;\; (x,p,y)\in V\},\;\;\mbox{for}\;(y,p)\in Y\times P.$$
For each $(y,p)\in Y\times P$, the lower semicontinuous envelope relative to $V$ of the function: $x\mapsto d(y,
F(x,p))$ is defined by
\begin{equation}
\varphi_V(x,y,p):=\left\{\begin{array}{lll}&\liminf_{u\to x,u\in V_{(y,p)}}d(y, F(u,p))&\quad\mbox{if}\quad x\in\cl V_{(y,p)}\\
&+\infty&\quad\mbox{otherwise.}
\end{array}\right.
\end{equation}
\vskip 0.2cm
Similarly to Theorem \ref{Char-strong-slope}, one has
\begin{theorem}\label{Char-strong-slope-bis}
Let $X$ be a complete metric space, $Y$ be a metric space and $P$ be a topological space. Let
  $F:X\times P\rightrightarrows Y$  be  a set-valued mapping and let $((x_0,p_0),y_0)\in \gph F;$ $V\subseteq X\times P\times Y;$  $\tau\in (0,+\infty)$  be given.  Suppose that
 for any $p$ near $\bar{p},$ the set-valued mapping $x \rightrightarrows
F(x,p)$ is a closed multifunction.
Then, among the following statements, one has $(i)\Leftrightarrow (ii)\Leftarrow (iii)\Rightarrow (iv).$ Moreover, if $Y$ is a normed linear space; $V_{(x,p)}$ is convex for any $(x,p)$ near $(x_0,p_0)$ and $V_{(y,p)}$ is open for $(y,p)$ near $(y_0,p_0),$ then $(i)\Rightarrow (iv).$

\item{(i)}  $F$ is metrically regular relative to $V$ uniformly in $p$ at $((x_0,p_0),y_0);$

\item{(ii)} There exist $\delta,\gamma>0$ and a neighborhood $W$ of $p_0$ such that
\begin{align}\label{global-char?}
|\Gamma\varphi_V(\cdot,y,p)|(x)\ge\tau^{-1}\;\;\mbox{for all}\;  (x,p,y)\in& \big(B(x_0,\delta)\times W\times B(y_0,\delta)\big),  \;x\in\cl V(y,p) \nonumber
\\&\;\mbox{with}\;d(y,F(x,p))\in (0,\gamma);
\end{align}
\item{(iii)} There exist $\delta,\gamma>0$ and a neighborhood $W$ of $p_0$ such that
\begin{align}\label{local-char@}
|\nabla\varphi_V(\cdot,y,p)|(x)\ge\tau^{-1}\;\;\mbox{for all}\;  (x,p,y)\in& \big(B(x_0,\delta)\times W\times B(y_0,\delta)\big),\;x\in\cl V(y,p)\nonumber\\ &\mbox{with}\;d(y,F(x,p))\in (0,\gamma);
\end{align}
\item{(iv)} There exist $\delta,\gamma>0$ and a neighborhood $W$ of $p_0$ such that
\begin{align}\label{local-char@@}
|\nabla\varphi_V(\cdot,y,p)|(x)\ge\tau^{-1}\;\;\mbox{for all}\;  (x,p,y)\in& \big(B(x_0,\delta)\times W\times B(y_0,\delta)\big)\cap V \nonumber\\&\mbox{with}\;d(y,F(x,p))\in (0,\gamma).
\end{align}
\end{theorem}
\vskip 0.2cm
\proof{Proof.}  The proof being similar  to the one of Theorem \ref{Char-strong-slope}, we omit it.\hfill\Halmos 
\section{Coderivative characterizations  of  directional metric regularity}
For the usual metric regularity, sufficient conditions in terms of   coderivatives have been given by various authors, for instance, in \cite{Aze06, JT1, Borisbook1, RefNT3}. In this section, we establish a characterization of  directional metric regularity  using the Fr\'{e}chet subdifferential in Asplund spaces,  i.e., Banach spaces for which every convex continuous function is generically Fr\'{e}chet differentiable. There are  many equivalent descriptions of Asplund spaces, which can be found, e.g., in \cite{Borisbook1}
and its bibliography. In particular, any reflexive space is Asplund, as well as each Banach space such that if each of its separable subspaces has a separable dual. 
%
%

In order to formulate 
in this section some coderivative characterizations of directional metric regularity, we require some more
definitions.
Let $X$ be  a Banach space. Consider now an extended-real-valued function  $f:X\rightarrow\R\cup\{+\infty\}.$  The {\it Fr\'{e}chet subdifferential} of $f$ at ${\bar x}\in\dom f$
is given as
\begin{eqnarray*}
 \partial f(\bar{x}) = \left\{ x^\ast \in X^\ast: \liminf_{x \rightarrow \bar{x}, \quad x \neq \bar{x}} \frac{f(x) -
f(\bar{x}) - \langle  x^\ast, x - \bar{x} \rangle}{\| x - \bar{x} \|} \ge
0 \right\}.
\end{eqnarray*}
 For convenience of the reader, we  would like to mention that  the terminology \textit{regular subdifferential}  instead of Fr\'echet subdifferential is also popular   due to its use in Rockafellar and Wets \cite{RockBook2002}.
Every element of the Fr\'echet subdifferential  is termed as a Fr\'echet (regular)  
subgradient. If $\bar{x} $ is a point where $ f( \bar{x}) = \infty $, 
then we set $ \partial f( \bar{x}) = \emptyset $. In fact one can show that an element $  x^\ast $ is a Fr\'echet  subgradient
of $ f $ at $ \bar{x} $ iff
$$
f(x) \ge f(\bar{x}) + \langle  x^\ast, x - \bar{x}\rangle + o (\| x - \bar{x}\|)\quad 
\text{where}\quad \displaystyle \lim_{x \rightarrow \bar{x}} \frac{o(\| x - \bar{x}\|)}{\|x - \bar{x}\|}
= 0.$$ 
It is well-known that the Fr\'{e}chet subdifferential satisfies a fuzzy sum
rule on Asplund spaces (see \cite[Theorem 2.33]{Borisbook1}). More precisely, if
$X$ is an Asplund space and $f_{1},f_{2}:X\rightarrow
\mathbb{R\cup\{\infty\}}$ are such that $f_{1}$ is Lipschitz continuous
around $\overline{x}\in\dom f_{1}\cap\dom %
f_{2}$ and $f_{2}$ is lower semicontinuous around $\overline{x},$
then for any $\gamma>0$ one has%
\begin{equation}
\partial(f_{1}+f_{2})(\overline{x})\subset%
{\displaystyle\bigcup}
\{\partial f_{1}(x_{1})+\partial f_{2}%
(x_{2})\mid x_{i}\in\overline{x}+\gamma\overline{B}_{X},\left\vert f
_{i}(x_{i})-f_{i}(\overline{x})\right\vert \leq\gamma,i=1,2\}+\gamma
B_{X^{\ast}}.\label{fuz}%
\end{equation}
For a nonempty closed set $C\subseteq X,$ denote by  $\delta_{C}$ the \textit{indicator function} associatedwith C
 (i.e. $\delta_{C}(x)=0$,  when $x\in C$ and  $\delta_{C}(x)=\infty$
otherwise). The {\textit{  Fr\'{e}chet normal cone} to $C$  at $\bar x$  is denoted by $N(C,\bar x)$. It  is a closed and convex object in $X^\ast$ which  is defined as
$\partial \delta_C(\bar x).$  Equivalently  a vector $ x^\ast \in X^\ast $ is a Fr\'echet normal to $ C $ at $ \bar{x}$ if
\begin{eqnarray*} \langle x^\ast , x - \bar{x}
\rangle \le o ( \| x - \bar{x} \|), \quad \forall x \in C,
\end{eqnarray*}
where $ \lim_{ x \rightarrow \bar{x}} \displaystyle\frac{o( \| x -
\bar{x} \|)}{\| x - \bar{x}\|} = 0 $. 
Let $F:X\rightrightarrows Y$ be a set-valued map and $(x,y)\in\gph F.$ Then the {\it Fr\'{e}chet coderivative} at
$(x,y)$ is the set-valued map $D^{\ast
}F(x,y):Y^{\ast}\rightrightarrows X^{\ast}$ given by
\[
D^{\ast
}F(x,y)(y^{\ast}):=\big\{x^{\ast}\in
X^{\ast}\mid(x^{\ast},-y^{\ast})\in N(\gph%
F,(x,y))\big\}.
\]
This  notion  is recognized as a powerful tool of variational analysis when applied to problems of optimization and control (see \cite{Borisbook1, Mo-Ou2007, kummer2002}, and the references therein).

In the proof of the main result, we will use the following particular version of 
 Theorem \ref{Char-strong-slope} for  directional metric regularity.
\vskip 0.5cm 
\begin{theorem}\label{Direct-Char-strong-slope}
Let $X$ be a complete metric space and $Y$ be a normed space. Let
  $F:X\rightrightarrows Y$  be  a closed multifunction (i.e., its graph is closed) and fix $(x_0,y_0)\in \gph F$ and   $ V\subseteq X\times Y.$  For  a given $\tau\in (0,+\infty),$ then among the following statements, one has $(i)\Leftrightarrow(ii)\Leftarrow (iii).$ 

\item{(i)}  $F$ is metrically regular in  the direction $\bar y$ at $(x_0,y_0);$

\item{(ii)} There exists $\delta>0$ such that
\begin{align}
|\Gamma\varphi_{V(\bar y,\delta)}(\cdot,y)|(x)\ge\tau^{-1}\;\;\mbox{for all}\;  (x,y)\in \big(B(x_0,\delta)&\times B(y_0,\delta)\big),\nonumber\\&x\in \cl V(\bar y,\delta)\;\; \mbox{with}\;d(y,F(x))\in (0,\delta);
\end{align}

\item{(iii)} There exist $\delta>0$ such that
\begin{equation}
|\nabla\varphi_{V(\bar y,\delta)}(\cdot,y)|(x)\ge\tau^{-1}\;\;\mbox{for all}\;  (x,y)\in \big(B(x_0,\delta)\times B(y_0,\delta)\big)\; x\in \cl V(\bar y,\delta) \;\;\mbox{with}\;d(y,F(x))\in (0,\delta).
\end{equation}

\end{theorem}
Denote by $S_{Y^*}$ the  unit  sphere in the dual space $Y^*$ of $Y,$ and by $d_*$ the metric associated with the dual norm on $X^*.$ For given $\bar y\in Y$ and $\delta>0,$ denote by
\begin{equation} \label{Dual}
C_{Y^*}(\bar y,\delta):= \{y^*\in Y^*:\;\; |\langle y^*,\bar y\rangle|\le\delta \}\; \mbox{and}\;S_{Y^*}(\bar y,\delta):= \{y^*\in Y^*:\;\; \|y^*\|\le 1+\delta,\; \langle y^*,\bar y\rangle\le\delta\},\end{equation}
and 
\begin{equation}
T(\bar y,\delta):=\{(y_1^*,y_2^*)\in S_{Y^*}(\bar y,\delta)\times C_{Y^*}(\bar y,\delta):\;\; \|y_1^*+y_2^*\|=1\}.
\end{equation}
\vskip 0.2cm
For a given multifunction $F: X\rightrightarrows Y,$ we associate the  multifunction  $G:X\rightrightarrows Y\times Y$ defined by
$$ G(x)=F(x)\times F(x),\quad x\in X. $$ 
Recall also that a multifunction $F: X\rightrightarrows Y$ is said to be \textit{ pseudo-Lipschitz} ( or  \textit{Lipschitz-like} or satisfying the Aubin property) around $(x_0,y_0)\in\gph F$ if there exist  constants $L,\delta>0$ such that
$$ F(x)\cap B(y_0,\delta)\subseteq F(x^\prime) +L\|x-x^\prime\| B_Y,\quad \mbox{for all}\;\; x,x^\prime\in B(x_0,\delta).$$
It is well known that $F$ is  pseudo-Lipschitz around $(x_0,y_0)$  if  and only if the function $d(\cdot,F(\cdot)): X\times Y\to \R$ is Lipschitz near $(x_0,y_0), \  $(see for instance \cite[Theorem 1.142]{penot2013}).
\vskip 0.2cm
A coderivative characterization of  directional metric regularity is initiated in the following theorem, which is  the main result of this section.
\begin{theorem}\label{Coder-Char}  Let $X,Y$ be Asplund spaces. Let $F:X\rightrightarrows Y$ be a closed multifunction  and $(x_0,y_0)\in \gph F$ be given. Let $F$ be pseudo-Lipschitz around $(x_0,y_0).$ Suppose that $F$ has convex values around $x_0$, i.e., $F(x)$ is convex for all $x$ near $x_0$. If 
\begin{equation}\label{Coder-DMR}
\liminf_{\substack{(x,y_1,y_2)\overset {G} {\rightarrow}(x_0,y_0,y_0)\\\delta\downarrow  0^+}}d_*(0,D^\ast G(x,y_1,y_2)(T(\bar y,\delta)))> m>0, 
\end{equation}
 then $F$ is directionally metrically regular in  the direction $\bar y$ with modulus $\tau\le m^{-1}$ at $(x_0,y_0).$The notation
$(x,y_1,y_2)\overset {G} {\rightarrow} (x_0,y_0,y_0)$ means that
$(x,y_1,y_2)\rightarrow(x_0,y_0,y_0)$ with $(x,y_1,y_2)\in
\gph G.$
\end{theorem}
\vskip 0.2cm  The following lemmata are  needed in the proof of Theorem \ref{Coder-Char}.
\begin{lemma}\label {localized-direct} Let $F:X\rightrightarrows Y$ be a multifunction with convex values for $x$ near $x_0$ and $(x_0,y_0)\in \gph F.$ Then for any $\bar y\in Y,$ $\delta_1,\delta_2>0,$ there exist $\eta,\delta>0$ such that for all $ x\in B(x_0,\eta),$ one has
\begin{equation}\label{Inclusion}
\big(F(x)+\cone B(\bar y,\delta)\big)\cap B(y_0,\eta)\cap \{y\in Y: \;d(y,F(x))<\eta\}\subseteq F(x)\cap B(y_0,\delta_1)+\cone B(\bar y,\delta_2).
\end{equation}
\end{lemma}
\vskip 0.2cm
\proof{Proof.}  Let $\bar y\in Y,$ $\delta_1,\delta_2$ be given. If $\|\bar y\|<\delta_2$ then the conclusion holds trivially. Suppose $\|\bar y\|\ge \delta_2.$ Take $\delta=\delta_2/2$ and $\varepsilon\in (0,\delta_1/2)$  sufficiently small such that
$$ \frac{\varepsilon(\|\bar y\|+\delta_2/2)}{\delta_1-2\varepsilon} <\delta_2/2.$$
Let $\eta\in (0,\varepsilon/2)$ such that $F(x)$ is convex for all $x\in B(x_0,\eta)$ and let
now $x\in B(x_0,\eta)$ and $y\in \big(F(x)+\cone B(\bar y,\delta))\cap B(y_0,\eta)$ with $d(y, F(x)) <\varepsilon/2$ be given. Then, there exist $z,v\in F(x)$ such that 
$$ y=z+\lambda (\bar y+\delta u),\;\;\mbox{for}\;\lambda\ge 0,\;u\in Y\;\;\|u\|\le 1;\;\;\|y-v\|<\varepsilon/2.$$
If $z\in B(y_0,\delta_1)$ then  the proof is over. Otherwise, one has
$$ \lambda(\|\bar y\|+\delta) \ge \|y-z\|\ge\|z-y_0\|-\|y-y_0\|\ge \delta_1-\eta>\delta_1-\varepsilon.$$ 
By setting
$$ t:=\frac{\delta_1-2\varepsilon}{\delta_1-\varepsilon},\;\; w:=tz+(1-t)v\in F(x),$$
one has
$$ \|w-y_0\| \le t\|z-v\| +\|v-y_0\| \le t\lambda\|y+\delta u\| +t\|y-v\| +\|v-y_0\| <\delta_1-2\varepsilon/2+\varepsilon<\delta_1$$
and,
$$ \frac{(1-t)\|y-v\| } {t\lambda}<\frac{\varepsilon(\|\bar y\|+\delta)}{\delta_1-2\varepsilon}<\delta_2/2.$$
Thus,
$$ y-w =t\lambda\left(\bar y+\delta u+(1-t)\frac{y-v}{t\lambda}\right)\in \cone B(\bar y,\delta_2),$$
which implies   that $y\in F(x)\cap B(y_0,\delta_1)+\cone B(\bar y,\delta_2).$\hfill\Halmos 
\vskip 0.5cm
Associated with the multifunction $F,$ for given $\varepsilon>0,$ $(x_0,y_0)\in\gph F,$ we define the \textit{localization}  of $F$ by
\begin{equation}
F_{(x_0,y_0,\varepsilon)}(x):=\:\left\{\begin{array}{ll}F(x)\cap \bar B(y_0,\delta_0)\quad&\mbox{if}\quad x\in \bar B(x_0,\varepsilon)\\
\emptyset\quad&\mbox{otherwise.}
\end{array}
\right.
\end{equation}
Note that, by definition, one has
\begin{equation}\label{Lo-co}
D^\ast F(x,y)=D^\ast F_{(x_0,y_0,\varepsilon)}(x,y)\quad\forall (x,y)\in \gph F\cap (B(x_0,\varepsilon)\times B(y_0,\varepsilon)).
\end{equation}
The preceding lemma implies  obviously  the next  corollary.
\begin{corollary}\label{Localization}
Let $F:X\rightrightarrows Y$ be a multifunction with convex values for $x$ near $x_0$ and $(x_0,y_0)\in \gph F.$ 
Then the two following conditions  are equivalent:
\begin{enumerate}
\item   $F$ is directionally metrically regular in  the direction $\bar y$;
\item  For any $\varepsilon>0, \, F_{(x_0,y_0,\varepsilon)}$  is directionally metrically  in the  direction $\bar y$.
\end{enumerate}
\end{corollary}
\begin{lemma} \label{Nor-cone}
Let $C\subseteq Y$ be a nonempty convex cone. One has
$$N(C,z)\subseteq \{z^*\in Y^*: \langle z^*,z\rangle= 0\},\quad\mbox{for all}\;\; z\in C.$$
\end{lemma}

\proof{Proof.}  Let $z\in C$ be given. Then $\lambda z\in C$ for all $\lambda>0.$ Hence, for $z^*\in N(C,z),$ one has
$$ \langle z^*,\lambda z-z\rangle\le 0\quad\forall \lambda>0.$$
Thus,  $\langle z^*,z\rangle=0.$  \hfill\Halmos 
\vskip 2mm
Given a  multifunction $F:X\rightrightarrows Y$, $\bar y\in Y,$  for $y\in Y$  and $\delta>0,$  we define the set $\mathcal{V}(y,\delta)\subseteq X\times Y$ by
\begin{equation} \label{def-V}
\mathcal{V}(y,\delta):=\{(x,z)\in X\times Y:\;\; y\in F(x)+z,\; z\in\cone \bar B(\bar y,\delta)\}.
\end{equation}
\begin{lemma}\label{Reg-V} Let $X,Y$ be Asplund spaces and let $F:X\rightrightarrows Y$ be a closed multifunction with convex values for $x$ near $x_0,$ and $(x_0,y_0)\in \gph F,$ $\bar y\in Y$ be given.   We suppose by assumption that 
\begin{equation}\label{Co-Re-V}
\lim_{\substack{(x,y)\overset  {F}\rightarrow (x_0,y_0)\\ \delta\downarrow  0^+}} d_*(0, D^\ast F(x,y)(C_{Y^*}(\bar y,\delta)\cap S_{Y^*}))>0.
\end{equation}
Then there exist $\kappa>0, \varepsilon >0, \delta_0>0$ such that for all $\delta\in (0,\delta_0),$ for any  $(x,y)\in B((x_0,y_0),\varepsilon), z\in \cone \bar B(\bar y,\delta),$ with $d(y,F(x))\in (0,\varepsilon_0) $ and $y-z\in F(x)\cap\bar B(y_0,\varepsilon),$  we can find $\eta>0$ such that
\begin{equation} d((x^\prime,z^\prime), \mathcal{V}(y,\delta)) \le \tau d(y, F(x^\prime)+z^\prime)\;\;\mbox{for all}\;\; (x^\prime,z^\prime)\in B((x,z),\eta),\; z^\prime\in\cone\bar B(\bar y,\delta).
\end{equation}
\end{lemma}

\textit{Proof.}
 Since (\ref{Co-Re-V}), then there exists $\delta_0\in (0,1)$ such that
\begin{equation} \label{Infer-Co}
\inf_{(x,y)\in\gph F\cap B((x_0,y_0),\delta_0)} d_*(0, D^\ast F(x,y)(C_{Y^*}(\bar y,\delta_0)\cap S_{Y^*})):=m>0.
\end{equation}
}
Let $\varPhi: X\times Y\rightrightarrows Y$ be a multifunction defined by 
$$ \varPhi(x,z):=\left\{ \begin{array}{ll} F(x) +z\quad&\mbox{if}\quad z\in \cone\bar B(\bar y,\delta),\\
\emptyset\quad\mbox{otherwise.}
\end{array}\right.$$
Then $\varPhi$ is a closed multifunction, and by a direct calculation, for $(x,z,y)\in\gph \varPhi,$ one has
\begin{equation}\label{Coder-Phi}
D^\ast \varPhi(x,z,y)(y^*)=\left\{(x^*,y^*+z^*)\in X^*\times Y^*: \;x^*\in D^\ast F(x,y-z)(y^*),\; z^*\in N(\cone \bar B(\bar y,\delta),z) \right\}.
\end{equation}
 Then, $\mathcal{V}(y,\delta)=\varPhi^{-1}(y).$ Let $\varepsilon,\delta\in (0,\delta_0/2)$ with $(1+\delta+|\bar y|)\delta<\delta_0.$
Let $(x,y)\in B((x_0,y_0),\varepsilon),$ $z\in \cone \bar B(\bar y,\delta)$ with $d(y,F(x))\in (0,\varepsilon);$ $y-z\in F(x)\cap\bar B(y_0,\varepsilon)$ be given. Take $\eta=\min\{\delta,\|z\|\}>0,$ and  $(x^\prime,z^\prime,y^\prime)\in \gph\varPhi \cap B((x,z,y),\eta),$  $(x^*,w^*)\in D^\ast \varPhi(x^\prime,z^\prime)(y^*)$ with $\|y^*\|=1.$  Then $x^*\in D^\ast F(x^\prime,y^\prime-z^\prime)(y^*);$ $w^*=y^*+z^*$ with some $z^*\in N(\cone \bar B(\bar y,\delta),z^\prime).$   Since $N(\cone \bar B(\bar y,\delta),z^\prime)\subseteq \{z^*\in Y^*: \langle z^*,z^\prime\rangle= 0\},$ then $|\langle z^*,\bar y\rangle|\le\delta\|z^*\|.$ If $\|w^*\|<\delta,$ then $\|z^*\|<1+\delta,$ 
and moreover,
 $$|\langle y^*,\bar y\rangle|\le |\langle z^*,\bar y\rangle| +\delta\|\bar y\|\le (1+\delta+\|\bar y\|)\delta<\delta_0.$$
Hence, $y^*\in C_{Y^*}(\bar y,\delta_0)\cap S_{Y^*}.$ Since $(x,y-z)\in B((x_0,y_0),\varepsilon),$ $(x^\prime,z^\prime,y^\prime)\in \gph\varPhi \cap B((x,z,y),\eta),$ then $(x^\prime,y^\prime-z^\prime)\in B((x_0,y_0),\delta_0)$
Therefore, from (\ref{Infer-Co}), we obtain $\|x^*\|\ge m$.
 Therefore,
$$ \liminf_{(x^\prime,z^\prime,y^\prime)\rightarrow_\varPhi(x,z,y)}d_*(0, D^\ast \varPhi(x^\prime,z^\prime,y^\prime)(S_{Y^*}))\ge \min\{m,\delta\}.$$
Thanks  to the standard coderivative characterization of   metric regularity for closed multifunctions (see, e.g., \cite{Aze06, Io00}), we conclude that $\varPhi$ is metrically regular around $(x,z,y).$  Thus,  there exists $\eta>0$ such that
$$ d((x^\prime,z^\prime), \mathcal{V}(y,\delta)) \le \tau d(y,\varPhi (x^\prime,z^\prime))= \tau d(y, F(x^\prime)+z^\prime)\;\;\mbox{for all}\;\; (x^\prime,z^\prime)\in B((x,z),\eta),\; z^\prime\in\cone\bar B(\bar y,\delta). $$
So the lemma is proved
.\hfill{$\Box$
\vskip 2mm
The next lemma is a penalty result which is similar to the one by Clarke (\cite{Clarke83}).
\begin{lemma}\label{penalty} Let $C$ be a subset of a metric space $X$ and let $x_0\in C$ and $\varepsilon>0.$ Then for a function $f:X\to\R\cup\{+\infty\}$ which is Lipschitz on $B(x_0,2\varepsilon)$ with constant $L>0,$ one has
$$f^*:= \inf\big\{f(x)\ :\; x\in C\cap B(x_0,2\varepsilon)\}\le  \inf\big\{f(x)+td(x,C)\ :\; x\in B(x_0,\varepsilon)\big\},$$
whenever $t\ge L.$
\end{lemma}
\proof{Proof.}  For any $x\in B(x_0,\varepsilon),$ pick a sequence $\{z_n\}_{n\in\N}\subseteq C$ with $\lim_{k\to\infty} d(x,z_k)=d(x,C).$ Then $z_k\in B(x_0,2\varepsilon)$ when $k$ is sufficiently large. Therefore,
$$ f^\ast \le f(z_k)\le f(x) +Ld(x,z_k)\rightarrow f(x)+Ld(x,C).\qquad\hfill \Box $$

\textit{Proof of Theorem \ref{Coder-Char}.} By the assumption, there is $\delta_0\in (0,1)$ such that
\begin{equation}\label{1bis}
\inf_{(x,y_1,y_2)\in\gph G\cap B((x_0,y_0,y_0), 2\delta_0)} d_*(0,D^\ast G(x,y_1,y_2)(T(\bar y,\delta_0)))\ge m+\delta_0.
\end{equation}
According to Corollary \ref{Localization} and relation (\ref{Lo-co}), by considering the localization $F_{(x_0,y_0,\delta_0)}$ instead of $F,$ without any  loss of generality, we can assume that
$$ F(x) \subseteq \bar B(y_0,\delta_0)\quad\mbox{for all}\;\; x\in \bar B(x_0,\delta_0).$$
Note that for all $(x,y_1)\in\gph F,$ $y_1^*\in Y^*,$ one has
$$  D^\ast G(x,y_1,y_1)((y^*_1,0))=D^\ast F(x,y_1)(T(y_1^*)).$$
Hence,
(\ref{Coder-DMR}) implies  obviously  (\ref{Co-Re-V}). Therefore, according to Lemma \ref{Reg-V}, there is $\kappa>0$ such that for all $\delta\in (0,\delta_0),$ for any $(x,y)\in B((x_0,y_0),\delta_0), \  z\in \cone \bar B(\bar y,\delta),$ with $d(y,F(x))\in (0,\delta_0)$ and $y-z\in F(x)$ (we may choose the     same $\delta_0$ as above), we can find $\gamma\in (0,\delta_0/2)$ such that
\begin{equation} \label{mathcal}d((x^\prime,z^\prime), \mathcal{V}(y,\delta)) \le \tau d(y, F(x^\prime)+z^\prime)\;\;\mbox{for all}\;\; (x^\prime,z^\prime)\in B((x,z),\gamma),\; z^\prime\in\cone\bar B(\bar y,\delta),
\end{equation}
where $ \mathcal{V}(y,\delta)$ is defined by (\ref{def-V}).
Since $F$ is pseudo-Lipschitz around $x_0,$ there is $\delta_0>0$  (we can assume it is the same than the previous one) and $L>0$ such that
\begin{equation}\label{PL}
F(x^\prime)\cap \bar B(y_0,\delta_0)\subseteq F(x) +L\|x-x^\prime\| B_Y\quad\forall x,x^\prime\in B((x_0,y_0),\delta_0).
\end{equation}
Moreover, the function $d(\cdot,F(\cdot)):X\times Y\to\R$  is Lipschitz around $(x_0,y_0)$  as recalled above,  say, on $B((x_0,y_0),\delta_0)$, with a Lipschitz  modulus  equal to $L.$ 
By  virtue of Theorem \ref{Direct-Char-strong-slope}, it suffices to show that  one has $|\nabla\varphi_\delta(\cdot,y)|(x)> m$  for any $(x,y)\in \big(B(x_0,\delta)\times B(y_0,\delta)\big)\; x\in \cl V_y(\bar y,\delta) \;\;\mbox{with}\;d(y,F(x))\in (0,\delta).$ Remind that, $V(\bar y,\delta),$ $V_y(\bar y,\delta),$ $\varphi_{\delta}(\cdot,y)$ are defined by (\ref{V}), (\ref{envelope}), respectively. Indeed, let $(x,y)\in B(x_0,\delta)\times B(y_0,\delta),  \ x\in \cl V(\bar y,\delta) \;\;\mbox{with}\;d(y,F(x))\in (0,\delta)$ be given. Set $|\nabla \varphi_\delta(\cdot,y)|(x):=\alpha.$ Since $d(y,F(\cdot))$ is Lipschitz on $B(x_0,\delta_0),$ then 
$$\varphi_\delta(x^\prime,y)=d(y,F(x^\prime))\quad\forall x^\prime\in B(x_0,\delta_0)\cap \cl V_y(\bar y,\delta).$$
By the definition of the
strong slope,  for each $\varepsilon\in (0,\delta),$ there is
$\eta\in (0,\varepsilon)$ with $$2\eta+\varepsilon<\min\{\gamma/2, \varepsilon d(y,F(x))\} \;\text{and}\;1-(\alpha+\varepsilon+2)\eta >0$$ such
that $$d(y,F(x^\prime))\geq (1-\varepsilon)d(y,F(x))\;\forall x^\prime\in
B(x,4\eta)$$ and that
$$d(y,F(x))\leq d(y,F(x^\prime))+(m+\varepsilon)\|x^\prime-x\|\quad \mbox{for all}\; x^\prime\in \bar B(x, 3\eta)\cap\cl V_y(\bar y,\delta).$$
Take $u\in B(x,\eta^2/4)\cap V_y(\bar y,\delta),$  $v\in F(u)$  such that $\|y-v\|\leq
d(y,F(x))+\eta^2/4.$ Then,
$$
\|y-v\|\leq d(y,F(x^\prime))+(\alpha+\varepsilon)\|x^\prime-x\|+\eta^2/4\;\forall x^\prime\in \bar B(u,2\eta)\cap \cl V_y(\bar y,\delta).
$$
Consequently, 
\begin{equation}\label{Penal}
\|y-v\|\leq d(y,F(x^\prime))+(\alpha+\varepsilon)\|x^\prime-u\|+(\alpha+\varepsilon+1)\eta^2/4\;\forall (x^\prime,z^\prime)\in \big(\bar B(u,2\eta)\times Y\big)\cap \mathcal{V}(\bar y,\delta).
\end{equation}
Let $z\in \cone B(\bar y,\delta)$ such that $y-z\in F(u).$  Then, 
 $$ \|z\| \ge d(y,F(u))\ge (1-\varepsilon)d(y,F(x))>0.$$
 Hence, by virtue of  relation (\ref{PL}),  there exists a neighborhood of $u,$ say $B(u,2\eta)$ such that
\begin{equation}\label{cone@}
\begin{array}{ll}
y\in F(u)\cap B(y_0,\delta_0)+z &\subseteq F(u^\prime)+\|u-u^\prime\|B_Y+z\\
&\subseteq F(u^\prime)+\cone B(\bar y,\delta_0)\quad\mbox{for all}\quad u^\prime \in B(u,2\eta).
\end{array}
\end{equation}
Since the function $d(y,F(\cdot))$ is Lipschitz  on $B(x_0,\delta_0),$ then from relation (\ref{Penal}), according to Lemma \ref{penalty},  it follows that there is $t>0$ such that 
$$ \|y-v\|\leq d(y,F(x^\prime))+(\alpha+\varepsilon)\|x^\prime-u\|+ td((x^\prime,z^\prime),\mathcal{V}(\bar y,\delta))+(\alpha+\varepsilon+1)\eta^2/4\;\forall (x^\prime,z^\prime)\in \bar B(u,\eta)\times \bar B(z,\eta).$$
Moreover, by (\ref{mathcal}), one obtains
$$\begin{array}{ll}  \|y-v\|\leq &d(y,F(x^\prime))+(\alpha+\varepsilon)\|x^\prime-u\|+ t\kappa d(y,F(x^\prime)+z^\prime)+(\alpha+\varepsilon+1)\eta^2/4\\
&\mbox{for all} (x^\prime,z^\prime)\in \bar B(u,\eta)\times \bar B(z,\eta),\; z^\prime\in\cone \bar B(\bar y,\delta).\end{array}$$
Thus, setting $G(x):=\;F(x)\times F(x),\; x\in X,$ we derive
$$ \begin{array}{ll}  \|y-v\|\leq \|y-w_1\|+&(\alpha+\varepsilon)\|x^\prime-u\|+ t\kappa \|y-w_2-z^\prime\|+\\&+\delta_{\gph G}(x^\prime,w_1,w_2)+\delta_{\cone \bar B(\bar y,\delta)}(z^\prime)+(\alpha+\varepsilon+1)\eta^2/4\\
&\mbox{for all} (x^\prime,w_1,w_2, z^\prime)\in \bar B(u,\eta)\times Y\times Y\times \bar B(z,\eta),\; z^\prime\in\cone \bar B(\bar y,\delta).
\end{array}$$
Next,  applying the
Ekeland variational principle to the function
$$\begin{array}{ll}(x^\prime,w_1,w_2, z^\prime)\mapsto \psi(x^\prime,w_1,w_2, z^\prime):=&\|y-w_1\|+(\alpha+\varepsilon)\|x^\prime-u\|+ t\kappa \|y-w_2-z^\prime\|+\\&+\delta_{\gph G}(x^\prime,w_1,w_2)+\delta_{\cone \bar B(\bar y,\delta)}(z^\prime)\end{array}$$
on $\bar B(u,\eta)\times Y\times Y\times \bar B(z,\eta),$ we can select $(u_1,v_1,v_2,z_1)\in
(u,v,y-z, z)+\frac{\eta}{4} B_{X\times Y\times Y\times Y}$ with $(u_1,v_1,v_2)\in \gph G,$ $z_1\in \cone\bar B(\bar y,\delta)$
such that \begin{equation}\label{Estim 1}\|y-v_1\|+\tau\kappa\|y-v_2-z_1\|\leq\|y-v\| (\leq
d(y,F(x))+\eta^2/4);\end{equation} and 
$$\psi(u_1,v_1,v_2,z_1)\le \psi(x^\prime,w_1,w_2, z^\prime) +(\alpha+\varepsilon+1)\eta\|(x^\prime,w_1,w_2, z^\prime)-(u_1,v_1,v_2,z_1)\|$$
for all $(x^\prime,w_1,w_2, z^\prime) \in \bar B(u,\eta)\times Y\times Y\times \bar B(z,\eta).$ 
Thus, $$0\in\partial (\psi+(\alpha+\varepsilon+1)\eta\|\cdot-(u_1,v_1,v_2,z_1)\|)(u_1,v_1,v_2,z_1).$$
According to the fuzzy sum rule, 
we can find $$v_3\in B(v_1,\eta);\; v_4\in B(v_2,\eta);$$
$$(u_2,w_1,w_2)\in
B(u_1,\eta)\times B(v_1,\eta\times B(v_2,\eta)\cap\gph G;\;z_2,\;z_3\in B(z,\eta);$$
$$v_3^*\in\partial\|y-\cdot\|(v_3);\;(u_2^*,-w_1^*,-w_2^*)\in
N(\gph G,(u_2,w_1,w_2));$$
$$(v_4^*,z_3^*) \in t\kappa\|y-\cdot-\cdot\|(v_4,z_3);\;\; z_2^*\in N(\cone \bar B(\bar y,\delta),z_2),$$ satisfying
\begin{equation}\label{Estim 2}\begin{array}{ll}&\|v_3^*-w_1^*\|<(\alpha+\varepsilon+2)\eta;\;\; \|v_4^*-w_2^*\|<(m+\varepsilon+2)\eta;\\
&\|z_3^*+z_2^*\|<(m+\varepsilon+2)\eta;\;\;\|u_2^*\|\leq \alpha+\varepsilon +(\alpha+\varepsilon+2)\eta.\end{array}\end{equation}
Since $v_3^*\in\partial\|y-\cdot\|(v_3) $ (note that $\|y-v_3\|\geq
\|y-v\|-\|v_3-v\|\geq d(y,F(x))-\varepsilon-2\eta>0$), then
$\|v^*_3\|=1$ and $\langle v_3^*,v_3-y\rangle=\|y-v_3\|.$ 
Thus,  
$\|w_1^*\|\le 1+(\alpha+\varepsilon+2)\eta,$ and
the first relation of (\ref{Estim 2}) follows that
$$
\langle w_1^*,w_1-y\rangle \geq \langle v_3^*,v_3-y\rangle -(\alpha+\varepsilon+2)\eta\|v_3-y\|-2\eta=(1-(\alpha+\varepsilon+2)\eta)\|v_3-y\|-2\eta.$$
As $\eta\le \varepsilon d(y,F(x))\le \varepsilon d(y,F(u))/(1-\varepsilon)$ for all $u\in B(x,\eta),$ one obtains
\begin{equation}\label{Estim3}
\langle w_1^*,w_1-y\rangle \geq(1-\varepsilon_1)\|w_1-y\|,
\end{equation}
where
$$ \varepsilon_1:=(\alpha+\varepsilon+2)\eta-2(\alpha+\varepsilon+2)\eta\varepsilon(1-\varepsilon)^{-1}-2\varepsilon(1-\varepsilon)^{-1}.$$
On the other hand, since $F(u_2)$ is convex and $w_1^*\in -N(F(u_2),w_1),$ then by relation (\ref{cone}), there is $w_1^\prime\in F(u_2)$ such that $y-w_1^\prime\in \cone B(\bar y,\delta_0).$ Therefore
$$\langle w_1^*,y-w_1^\prime\rangle=\langle w_1^*,y-w_1^\prime\rangle+\langle w_1^*,w_1-w_1^\prime\rangle<0.$$
Consequently,
\begin{equation}\label{Estim-bis}
\langle w_1^*,\bar y\rangle \le \delta_0\|w_1^*\|/2\le \delta_0(1+(\alpha+\varepsilon+2)\eta)/2 .
\end{equation}
\vskip 0.2cm
Next, since $(v_4^*,z_3^*) \in t\kappa\|y-\cdot-\cdot\|(v_4,z_3),$ then $v_4^*=z_3^*$ and
$ \|z_3^*\|\le  t\kappa.$ 
Hence from (\ref{Estim 2}), one has
$$ \|w_2^*-z_2^*\|\le\|w_2^*-v_4^*\|+\|z_2^*-z_3^*\|\le 2(\alpha+\varepsilon+2)\eta.$$
As $z^*_2\in N(\cone\bar B(\bar y,\delta),z_2),$ with $z_2\not= 0$, then $\langle z_2^*,z_2\rangle=0.$
Therefore,
$$ |\langle w_2^*,z_2\rangle| \le 2(\alpha+\varepsilon+2)\eta\|z_2\|<.$$  As $z_2\in \cone \bar B(\bar y,\delta),$ one obtains
\begin{equation}\label{Es}
|\langle w_2^*,\bar y\rangle|\le 2(\alpha+\varepsilon+2)\eta +\delta.
\end{equation}
Moreover,
$$ |\langle w_2^*,w_2-y\rangle| \le |\langle z_2^*,w_2-y-z_2\rangle| +|\langle z_2^*-w_2^*,w_2-y\rangle|\le \varepsilon_2\|w_1-y\| ,$$
where
 $$\varepsilon_2=\big((t\kappa+2(\alpha+\varepsilon+2)+2(\alpha+\varepsilon+2)(\|y_0\|+2\delta_0+2\eta)\big)\varepsilon(1-\varepsilon)^{-1}. $$
The second inequality of the preceding relation follows from
$$\|z_2^*\|\le\|z_3^*\| +\|z_2^*-z_3^*\|\le t\kappa+(\alpha+\varepsilon+2),$$
and
$$ \|w_2-y\|\le \|w_2-v_2\|+\|v_2-(y-z)\| +\|z\|<2\eta+2\delta_0+\|y_0\|.$$
Hence, by using the convexity of $F(u_2),$ and $w_2^*\in -N(F(u_2), w_2)$
\begin{equation}\label{Estim 4}
\langle w_2^*,w_1-y\rangle=\langle w_2^*,w_1-w_2\rangle+\langle w_2^*,w_2-y\rangle\ge -\varepsilon_2\|w_1-y\|.
\end{equation}
From relations (\ref{Estim 3}) and (\ref{Estim 4}), one derives that
\begin{equation}\label{Estim 5}
 \langle w_1^*+w_2^*,w_1-y\rangle \ge (1-\varepsilon_1-\varepsilon_2)\|w_1-y\|.
\end{equation}
Consequently, $\|w_1^*+w_2^*\|\ge 1-\varepsilon_1-\varepsilon_2.$
\vskip 0.2cm
Set 
 $$y_1^*=\frac{w^*_1}{\|w_1^*+w_2^*\|};\;\; y_2^*=\frac{w^*_2}{\|w_1^*+w_2^*\|}\;\;\mbox{and}\;\; x^*=\frac{u^*_2}{\|w_1^*+w_2^*\|}.$$  From relations (\ref{Estim-bis}), (\ref{Es}), (\ref{Estim 5}), one
has   
$$ \langle y_1^*,\bar y\rangle\le  \delta_0\|y_1^*\|\le \frac{\delta_0(1+(\alpha+\varepsilon+2)\eta)}{2(1-\varepsilon_1-\varepsilon_2)} ;$$
$$ \langle y_2^*,\bar y\rangle\le \frac{2(\alpha+\varepsilon+2)\eta +\delta}{1-\varepsilon_1-\varepsilon_2};$$
$$x^*\in D^\ast G(u_2,w_1,w_2)(y_1^*,y_2^*);\;\; \|y_1^*+y_2^*\|=1.$$ As $\varepsilon_1,\varepsilon_2$ go to $0$ as $\varepsilon,\eta\to 0,$ then $y_1^*\in S_{Y^*}(\bar y,\delta)$ and $y_2^*\in C_{Y^*}(\bar y,\delta).$ Thus $(y_1^*,y_2^*)\in T(\bar y,\delta).$ As $(u_2,w_1,w_2)\in B((x_0,y_0,y_0),\delta_0),$ according to  (\ref{Estim 2}), one obtains
\begin{equation}\label{Estim 3}
\alpha+\delta_0\le\|x^*\|=\|u^*_2\|/\|w_1^*+w_2^*\|\leq \frac{\alpha+\varepsilon
+(\alpha+\varepsilon+2)\eta}{1-\varepsilon_1-\varepsilon_2}.\end{equation} 
As $\varepsilon,\eta,\varepsilon_1,\varepsilon_2$ are arbitrary small, we obtain $m+\delta_0\le\alpha$ and  the proof is complete.\hfill$\Box$
\vskip 0.2cm
Condition (\ref{Coder-DMR}) is also a necessary condition for  directional metric regularity in Banach spaces as showed in the next  proposition.
\begin{proposition}\label{nec-cond} Let $X,Y$ be Banach spaces and let  $F:X\rightrightarrows Y$ be a closed multifunction, $(x_0,y_0)\in \gph F$  and $\bar y\in Y.$  Suppose that $F$ has convex values for $x$ near $x_0.$ If 
$F$ is metrically  regular in the direction $\bar y\in Y$ at $(x_0,y_0)$, then 
\begin{equation}\label{Coder-DMR-bis}
\liminf_{\substack{(x,y_1,y_2)\overset {G}\rightarrow(x_0,y_0,y_0)\\ \delta\downarrow 0^+}}d_*(0,D^\ast G(x,y_1,y_2)(T(\bar y,\delta)))>0. . \notag
\end{equation}
\end{proposition}
\vskip 0.2cm
\proof{Proof.}  Assume that $F$ is metrically  regular in  the direction $\bar y\in Y,$ i.e.,  there exist $\tau>0,\delta>0,\varepsilon>0$ such that
\begin{equation}\label{DR}
 d(x,F^{-1}(y)) \le \tau d(y, F(x))\quad\mbox{for all}\;\; (x,y)\in B(x_0,\varepsilon)\times B(y_0,\varepsilon);\; y\in F(x)+\cone B(\bar y,\delta).
\end{equation}
For $\gamma\in(0,\delta),$ 
let $(x,y_1,y_2)\in \mbox{gph}G\cap B(x_0,\varepsilon/2)\times
B(y_0,\varepsilon/2)\times B(y_0,\varepsilon/2);$ $(y_1^*,y^*_2)\in T(\bar y,\gamma)$ and $x^*\in D^\ast G(x,y_1,y_2)(y_1^*,y_2).$ For any
$\alpha\in (0,1),$ there exists $\beta\in (0,\varepsilon/2)$ such that
\begin{equation}\label{Coderi}\langle x^*,u-x\rangle-\langle
y_1^*,v_1-y_1\rangle+\langle
y_2^*,v_2-y_2\rangle\leq \varepsilon(\|u-x\|+\|v_1-y_1\|+\|v_2-y_2\|),\end{equation}
for all $(u,v_1,v_2)\in\mbox{gph}G\cap B((x,y_1,y_2),\beta).$
\vskip 0.2cm
  For $\delta_1\in (0,\delta),$ take $w\in B_{Y}$ such that $\langle y_2^*,\bar y+\delta w\rangle \le \gamma-\delta_1.$  Since (\ref{DR}), for all sufficiently small $t>0,$ we can find $u\in 
F^{-1}(y_2+t(\bar y+\delta w))$ such that
$$\|x-u\|\leq (1+\alpha)\tau d(y_2+t(\bar y+\delta w),F(x))\leq (1+\alpha)t\|\bar y+\delta u\|<\beta.$$
 Since $y_1^*\in -N(F(x),y_1)$ and $F(x)$ is convex, then $\langle y_1^*,y_2-y_1\rangle\ge 0.$ Therefore, by taking $v_1=v_2=y_2+t(\bar y+\delta w)$ into account in (\ref{Coderi}), one obtains
\begin{equation}\begin{array}{ll}
(1+\alpha)\tau t\|\bar y+\delta u\|\|x^*\|\geq\langle x^*,x-u\rangle &\geq \langle y_1^*+y_2^*,v-y_2\rangle+\langle y_1^*,y_2-y_1\rangle\\
&\geq t(\delta_1-\gamma)-\alpha t\|\bar y+\delta u\|((1+\alpha)\tau+1).\end{array}\notag
\end{equation}
As $\alpha>0,$ $\delta_1\in (0,\delta)$ are arbitrary, one has 
$$\|x^*\|\geq \frac{\delta-\gamma}{\tau\|\bar y+\delta u\|}\ge \frac{\delta-\gamma}{\tau(\|\bar y\|+\delta)}.
$$
Thus,
$$\liminf_{\substack{(x,y_1,y_2)\overset {G}\rightarrow (x_0,y_0,y_0)\\\gamma\to 0^+}}d_*(0, D^\ast G(x,y_1,y_2)(T(\bar y,\gamma)))\ge \frac{\delta}{\tau(\|\bar y\|+\delta)}>0. $$ The  proof is complete.
\hfill$\Box$
\vskip 0.2cm
Combining this proposition and Theorem \ref{Coder-Char}, one has
\begin{theorem}\label{Asplund}
Let $X,Y$ be Asplund spaces.  Suppose  $F:X\rightrightarrows Y$ be a closed multifunction and $(x_0,y_0)\in \gph F$ such that $F$ has convex values around $x_0.$ Suppose further that $F$ is pseudo-Lipschitz around $(x_0, y_0).$ Then, $F$ is metrically  regular in direction $\bar y\in Y$ at $(x_0,y_0)$ if and only if
\begin{equation}\label{Coder-DMR-bis1}
\liminf_{\substack{(x,y_1,y_2)\overset {G}\rightarrow(x_0,y_0,y_0)\\ \delta\downarrow 0^+}}d_*(0, D^\ast G(x,y_1,y_2)(T(\bar y,\delta)))>0. \notag
\end{equation}
\end{theorem}
Recall that the {\it Mordukhovich limiting coderative} of $F$ denoted by $D^\ast _MF(x,y):Y^*\rightrightarrows X^*$ is defined by
\begin{equation}D^\ast_MF(x,y)(y^\ast) :=\liminf_{\substack{u,v)\overset{F}\rightarrow (x,y)\\ 
v^*\overset{w^*}\rightarrow  y^*} }D^\ast F(u,v)(v^*)
=\left\{x^*\in X^*:\;\; \begin{array}{ll}&(x_n,y_n)\overset{F}\rightarrow (x,y)\\
&x^*_n\in  D^\ast F(x_n,y_n)(y_n^*)\\
&y_n^*\overset{w^*}\rightarrow  y^*,\; x_n^*\overset{w^*}\rightarrow  x^*\end{array} \right\}.\end{equation}

Let us now recall the notion of    \textit{partial sequential normal compactness}   (PSNC, in short, see \cite[page 76]{Borisbook1}). 
A  multifunction $F:X\rightrightarrows Y$  is \textit{partially sequentially normally compact}   at $(\bar{x},\bar{y})\in \gph F$, iff,  for any sequences  $\{(x_k , y_k , x_{k}^{\star} ,y_{k}^{\star})\}_{n\in \N}\subset  \gph F \times X^{\star} \times Y^{\star}$
satisfying $$(x_k , y_k) \to (\bar{x}, \bar{y}), x_{k}^{\star} \in D_M^{\star}F(x_k , y_k)(y_{k}^{\star}), x^{\star}_{k}\overset {w^{\star}} {\rightarrow} 0, \Vert y_{k}^{\star}\Vert\to 0,$$ 
one has $\Vert x^{\star}_{k}\Vert \to 0$ as $k\to \infty$.
\begin{remark}
Condition (PSNC) at $(\bar{x},\bar{y})\in \gph F$ is satisfied if $X$ is finite dimensional, or  $F$ is pseudo-Lipschitz around that point.
\end{remark} 
\vskip 2mm
The next  corollary that follows directly from the preceding theorem, gives a point-based condition for  directional metric regularity.
\begin{corollary}
Under the assumptions of Theorem \ref{Asplund}, suppose further that $G^{-1}$ is PSNC at $(x_0,y_0,y_0).$ Then $F$ is  metrically  regular in the direction $\bar y\in Y$ at $(x_0,y_0)$ if and only if
\begin{equation}
d_\ast(0,D_M^\ast G(x_0, y_0,y_0)(T(\bar y,0))) > 0 \nonumber. 
\end{equation}
\end{corollary}
With an analogous proof, we obtain the following parametric version  of  Theorem \ref{Coder-Char} . 
\begin{theorem}\label{Para-Coder} Let $X,Y$ be Asplund spaces and $P$ be a topological space. Let
  $F:X\times P\rightrightarrows Y$  be  a set-valued mapping and let $((x_0,p_0),y_0)\in \gph F.$ Suppose the following conditions are satisfied:

\item{(a)} For any $p$ near $\bar{p},$ the set-valued mapping $x \rightrightarrows
F(x,p)$ is a closed multifunction;

\item{(b)} For $(x,p)$ near $(x_0,p_0),$ $F(x,p)$ is convex;

\item (c) $F(\cdot,p)$ is pseudo-Lipschitz uniformly in $p$ around $(x_0,p_0).$
\vskip 0.2cm
\noindent Then, for $\bar y\in Y,$ $F$ is directionally metrically regular in direction $\bar y$ uniformly in $p$ at $(x_0,p_0,y_0)$  if and only if
\begin{equation}\label{Coder-DMR-para}
\liminf_{\substack{(x,p,y_1,y_2)\overset{G}\rightarrow (x_0,p_0,y_0,y_0)\\ \delta\downarrow  0^+}}d_*(0,D^\ast G_p(x,y_1,y_2)(T(\bar y,\delta)))>0, 
\end{equation}
where,
$$ G(x,p)=G_p(x):=F(x,p)\times F(x,p),\quad (x,p)\in X\times P,$$
\end{theorem}
We next consider a special case of $F(x,p):= f(x,p) -K:=f_p(x)-K,$ here, $K\subseteq Y$ is a nonempty closed convex subset. $f: X\times P\rightarrow Y$ is a (locally) continuous mapping around a given point $(x_0,p_0)\in X\times P$ with $f(x_0,p_0)\in K,$ and $f(\cdot,p)$ is Lipschitz uniformly in p near $(x_0,p_0).$ Obviously, for this case, assumptions $(a), (b), (c)$ of Theorem \ref{Para-Coder} are satisfied as well. Moreover, by setting $g_p:=(f_p,f_p):X\to Y\times Y,$ one has
$$ D^\ast G_p(x,y_1,y_2)(y^*)=\left\{\begin{array}{ll}D^\ast g_p(x)(y^*)\;&\mbox{if}\;\; f(x,p)-y_i\in K,\;y_i^*\in N(K, f(x,p)-y_i),\;i=1,2\\
\emptyset&\;\;\mbox{otherwise},
\end{array} \right.$$
where, we use the usual notations: $f_p(x):=f(x,p);$ $D^\ast f_p(x)(y^*):=D^\ast f_p(x, f(x,p))(y^*).$ Hence, Theorem \ref{Para-Coder} yields the following corollary.
\begin{corollary}\label{Coro} Let $X,Y$ be Asplund spaces and $P$ be a topological space. Let $K\subseteq Y$ be a nonempty closed convex subset and let
  $f:X\times P\rightarrow Y$ be  a locally continuous mapping around $(x_0,p_0)\in X\times P$ with $k_0:=f(x_0,p_0)\in Q.$  Suppose further that $f(\cdot,p)$ is Lipschitz uniformly in p near $(x_0,p_0).$
If for $\bar y\in Y,$ 
\begin{equation}\label{Convex-Cond}
\liminf_{\substack{(x,p,k_1,k_2)\rightarrow(x_0,p_0,k_0,k_0) \\\delta\downarrow  0^+}}d_*(0,D^\ast f_p(x)(T(\bar y,\delta))\cap  (N(K,k_1)\times N(K,k_2)))> m>0, 
\end{equation}
 then the mapping $F(x,p):=f(x,p)-K,\; (x,p)\in X\times P$ is directionally metrically regular in direction $\bar y$ uniformly in $p,$ with modulus $\tau=m^{-1}$ at $(x_0,p_0),$ i.e., there exist $\varepsilon>0,$ $\delta>0$ and a neighborhood $W$ of $p_0$such that
$$ d(x, S(y,p)) \le\tau d(f(x,p)-y, K)\quad\mbox{for all}\;\; (x,p,0)\in B(x_0,\varepsilon)\times W\times B(0,\varepsilon),$$
with $y\in f(x,p)-K+\cone B(\bar y,\delta).$ 
\vskip 0.2cm
 In particular, one has
$$ d(x, S(p)) \le\tau d(f(x,p), K)\quad\mbox{for all}\;\; (x,p)\in B(x_0,\varepsilon)\times W,$$
with $ f(x,p)\in K-\cone B(\bar y,\delta).$ Here, 
$$ S(y,p) =\{x\in X:\;\; f(x,p)-y\in K\},\quad S(p):=\{x\in X:\;\; f(x,p)\in K\}.$$
\end{corollary}
\begin{remark}\label{re}
Note that if $K$ is \textit{ sequentially normally compact at $\bar k$},  i.e., for all sequences $(k_n)_{n\in\N}\subseteq K,$ $(k^*_n)_{n\in\N}$ with $k^*_n\in N(K,k_n),$
$$ k_n\to \bar k \ \;\; k^*_n\overset{w^*}\rightarrow 0\quad \iff \|k^*_n\|\rightarrow 0,$$
 and $P$ is a \textbf{metric space}, then instead of (\ref{Convex-Cond}), the following point-based condition
\begin{equation}
d_*(0,D_{{\rm lim}}^*g_{p_0}(x_0)(T(\bar y,0)\cap (N(K,k_0)\times N(K,k_0)))>0, 
\end{equation}
is also a sufficient condition for  directional metric regularity at $\bar y,$ uniformly in $p$ of $F(x,p):=f(x,p)-K$ at $(x_0,p_0).$ Here, $D_{{\rm lim}}^*g_{p_0}(x_0)$ denotes the  sequential limiting subdifferential of $D^\ast  g_p(x):$
$$ D_{{\rm lim}}^*g_{p_0}(x_0)(y_1^*,y_2^*):=\liminf_{\substack{(x,p)\to (x_0,p_0)\\  (z_1^*,z_2^*)\overset{w^*}\rightarrow (y^*_1,y_2^*)}} D^\ast g_p(x)(z^*_1,z_2^*),\;\; y_1^*,y_2^*\in Y^*.$$
\end{remark}
\vskip 0.2cm
\begin{corollary}\label {Diff} With the assumptions of Corollary \ref{Coro}, suppose further that $f$ is Fr\'{e}chet differential  with respect to $x$ near $(x_0,p_0),$ and its derivative with respect to $x$ is continuous at $(x_0,p_0).$ Then, the mapping $F(x,p):=f(x,p)-K,\; (x,p)\in X\times P$ is directionally metrically regular in direction $\bar y$ uniformly in $p$ if and only if
\begin{equation}\label{Diff-Cond}
\liminf_{\substack{(k_1,k_2)\rightarrow k_0 \\ \delta\downarrow  0^+}}d_*(0,g^{\prime*}_x(x_0,p_0)(T(\bar y,\delta)\cap (N(K,k_1)\times N(K,k_2)))> m>0. 
\end{equation}
Here, $f^{\prime*}_x(x,p)$ stands for the adjoint operator of $f^{\prime}_x(x,p)$
Moreover, if $K$ is  normally  sequentially compact, then (\ref{Diff-Cond}) is equivalent to 
\begin{equation}\label{Limit-Diff-Cond}
d_*(0,g^{\prime*}_x(x_0,p_0)(T(\bar y,0)\cap (N(K,k_0)\times N(K,k_0)))>0. 
\end{equation}

\end{corollary}
\vskip 0.2cm
\proof{Proof.}  For the sufficiently part, suppose that 
$$ \liminf_{\substack{(k_1,k_2)\rightarrow (k_0,k_0) \\ \delta\downarrow 0^+}}d_*(0,g^{\prime*}_x(x_0,p_0)(T(\bar y,\delta))\cap N(K,k_1)\times N(K,k_2))> m>0. 
 $$ 
Since $f^\prime_x$ is continuous at $(x_0,p_0),$ for any $\varepsilon>0,$ there exist $\delta>0$ and a neighborhood $W$ of $p_0$ such that 
$$ \|g^\prime_x(x,p)-g^\prime_x(x_0,p_0)\| <\varepsilon\quad\mbox{for all}\;\; (x,p)\in B(x_0,\varepsilon)\times W.$$
Therefore, for all $\delta>0,$
$$ \|g^\prime_x(x,p)(y_1^*,y_2^*)-g^\prime_x(x_0,p_0)(y_1^*,y_2)\| <\varepsilon,$$
$$\mbox{for all}\;\; (x,p)\in B(x_0,\varepsilon)\times W,\; k_1,k_2\in B(k_0,\varepsilon),\;(y_1^*,y_2^*)\in T(\bar y,\delta)\cap (N(K,k_1)\times N(K,k_2)).$$
Consequently,
$$\begin{array}{ll}&\displaystyle \liminf_{\substack{(x, p,k_1, k_2)\rightarrow (x_0, p_0, k_0, k_0) \\  \delta\downarrow 0^+}}d_*(0,g^{\prime*}_x(x,p)(T(\bar y,\delta)\cap (N(K,k_1)\times N(K,k_2)))\\&=\displaystyle\liminf_{\substack{k\rightarrow k_0\\  \delta\downarrow  0^+}}d_*(0,g^{\prime*}_x(x_0,p_0)(T(\bar y,\delta)\cap (N(K,k_1)\times N(K,k_2)))> m>0 \end{array}.$$
The conclusion follows from Corollary \ref{Coro}. The proof of the necessary part is analogous to the one of Proposition \ref{nec-cond}. The equivalence between (\ref{Diff-Cond}) and (\ref{Limit-Diff-Cond}) follows from Remark \ref{re}.
\hfill\Halmos 
\vskip 0.2cm
Corollary \ref{Diff} subsumes  the  following result, established by Arutyunov, Avakov and Izmailov in \cite{ArAvIz07}.
\begin{corollary}\label{AAI}(\cite{ArAvIz07}, Theorem 2.3)
With the assumptions of Corollary \ref{Diff}, if
\begin{equation}\label{Rob-Type}
\cone\{\bar y\}\cap\int\left(f(x_0,p_0)+\im f^\prime(x_0,p_0)-K\right)\not=\emptyset,
\end{equation}
then the mapping $F(x,p):=f(x,p)-K,\; (x,p)\in X\times P,$ is directionally metrically regular in direction $\bar y$ uniformly in $p$ at $(x_0,p_0).$
\end{corollary}
\vskip 0.2cm
\proof{Proof.}  It suffices to show that (\ref{Rob-Type}) implies (\ref{Diff-Cond}). Indeed, assume (\ref{Rob-Type}) holds, and assume to contrary that (\ref{Diff-Cond}) fails to be hold. Then, there exist sequences $(\delta_n)_{n\in\N}$ with $\delta_n\downarrow 0;$ $(k^1_n)_{n\in\N}, (k_n^2)_{n\in\N}\subseteq K,\;k^i_n\to k_0=f(x_0,p_0)\; (i=1,2), \;(y^{1*}_n)_{n\in\N}, (y^{2*}_n)_{n\in\N}$ with $(y^{1*}_n,y^{2*}_n)\in T(\bar y,\delta_n)\cap [N(K,k^1_n)\times N(K,k^2_n)]$ and $(x_n^*)_{n\in \N}\subseteq X^*$ such that
$$ x_n^*=(y_n^{1*}+y_n^{2*})\circ f^\prime_x(x_0,p_0);\;\; \|x_n^*\| \to 0.$$
By (\ref{Rob-Type}), there exist $\lambda\ge 0,$  such that 
\begin{equation}\label{Inte}
0\in \int (f(x_0,p_0)+\im f^\prime(x_0,p_0)-K) -\lambda\bar y. 
\end{equation}
Set  $C_{nm}:= \;nf^\prime(x_0,p_0)(B(0,1)) +m(f(x_0,p_0)-K-\lambda\bar y),\;\;n\in\N.$
Then, $ \bigcup_{n,m=1}^\infty C_{nm}=Y.$
According to  the Baire   theorem, at least one of  the  $\cl C_{nm}'s$ has a nonempty interior.  Therefore,  consider  $y\in Y,  \; \alpha>0$ and $\varepsilon>0$ such that 
$$B(y,\varepsilon)\subseteq \cl(f^\prime(x_0,p_0)(B(0,\alpha)) +f(x_0,p_0)-K-\lambda\bar y).$$ 
On the  other hand, from (\ref{Inte}), there are $t, r>0$ such that 
$$-ty\in f^\prime(x_0,p_0)(B(0,\alpha)) +r(f(x_0,p_0)-K-\lambda\bar y).$$
Hence,
$$ B(0, t\varepsilon) \subseteq -ty+B(ty,t\varepsilon)\subseteq\cl((1+t)f^\prime(x_0,p_0)(B(0,\alpha)) +(t+r)(f(x_0,p_0)-K-\lambda\bar y).$$
Equivalently,  for $\gamma:=t\varepsilon/(t+r),$ $\beta:=(1+t)\alpha/(t+r),$
\begin{equation}\label{Interior}
B(\lambda\bar y, \gamma) \subseteq \cl(f^\prime(x_0,p_0)(B(0,\beta)) +f(x_0,p_0)-K).
\end{equation}
For each $n,$ let $u_n\in B_X$  be chosen such that $\langle y_n^{1*}+y_n^{2*},u_n\rangle<-1/2.$ Since $(y^{1*}_n,y^{2*}_n)\in T(\bar y,\delta_n),$ then
$$ \limsup_{n\to\infty} \langle y_n^{1*}+y_n^{2*},\lambda\bar y+\gamma u_n\rangle\le -\gamma/2.$$
On the other hand, by (\ref{Interior}), for each $n,$ we can find $x_n\in B(0,\beta),$ $z_n\in K$ such that
$$ \|\lambda\bar y+ \gamma u_n-(k_0+f_x^\prime(x_0,p_0)(x_n)-z_n)\|<1/n.$$
Since $y^{i*}_n\in N(K,k^i_n)\;(i=1,2);$ $\|x_n^*\|\to 0;$ $(x_n)$ is bounded, and $k^i_n\to k_0,$ one has
$$\limsup_{n\to\infty} \langle y_n^{1*}+y_n^{2*},\lambda\bar y+\gamma u_n\rangle= \limsup_{n\to\infty}(\langle y_n^{1*}+y_n^{2*},k_0-z_n\rangle +\langle x^*_n,x_n\rangle)\ge 0, $$
a contradiction. \hfill\Halmos 
\vskip 0.2cm
 From Theorem \ref{Coder-Char},  we can derive directly the following result due to  Ioffe (\cite{RefIoffe}) on  directional metric regularity of a closed convex multifunction for the case in which the convex multifunction under consideration  is assumed to be  pseudo-Lipschitz.

\begin{corollary}\label{Convex-Multi}(\cite{RefIoffe}, Proposition 15) Let $X,Y$ be  Banach spaces and $F:X\rightrightarrows Y$ be a closed convex multifunction and let $(x_0,y_0)\in \gph F$ and $\bar y\in Y.$ Suppose that $x_0\in \int F^{-1}(Y).$ Then $F$ is directionally metrically regular in direction $\bar y$ at $(x_0,y_0)$  if and only if 

\begin{equation}\label{cone}
\cone\{\bar y\}\cap\int(F(X)-y_0)\not=\emptyset.
\end{equation}
\end{corollary}
\vskip 0.2cm
\proof{Proof.}   For the sufficiency part, under the assumption $x_0\in \int F^{-1}(Y),$ according to the Robinson-Ursescu Theorem  (\cite{robinson73}, \cite{urcescu}), then $F$ is pseudo-Lipschitz near $x_0.$ By virtue of Theorem \ref{Coder-Char}, we only need to show that
\begin{equation}\label{CC}
\liminf_{\substack{(x,y_1,y_2)\overset{G}\rightarrow (x_0,y_0,y_0) \\ \delta\downarrow 0^+}}d_*(0,D^\ast G(x,y_1,y_2)(T(\bar y,\delta)))>0.
\end{equation}
By (\ref{cone}), there is  some $\lambda\ge 0$ such that $0\in \int (F(X)-y_0-\lambda\bar y).$ Then, by the convexity of the multifunction $F,$
$$ \bigcup_{n,m=1}^\infty (F(B(x_0,n))-m(y_0-\lambda\bar y))=Y.$$
Thanks to the Baire Category Theorem, similarly to the proof of Corollary \ref{AAI}, we can find $\gamma>0,$ $\beta>0$ such that
\begin{equation}\label{beta}
 B(\lambda\bar y,\gamma) \subseteq \cl(F(B(x_0,\beta))-y_0).
\end{equation}
Let $(x_n,y^1_n,y_n^2)\in \gph G;$ $(x^*_n)_{n\in\N},$ $(y^{1*}_n)_{n\in\N},$ $(y^{2*}_n)_{n\in\N},$ $(\delta_n)_{n\in\N}$ such that
$$(x_n,y^1_n,y_n^2)\to (x_0,y_0,y_0);\;\;\delta_n\downarrow 0^+;\;\; x^*_n\in D^\ast F(x_n,y^1_n,y^2_n)(y^{1*}_n,y^{2*}_n);\;\; (y^{1*}_n,y^{2*}_n)\in T(\bar y,\delta_n).$$
For each $n,$ take $u_n\in B_Y$ such that $\langle y^{1*}_n+y^{2*}_n,u_n\rangle <-1+\delta_n.$ Then, 
$$\langle y^{1*}_n+y^{2*}_n,\lambda\bar y+\gamma u_n\rangle<\delta_n-(1-\delta_n)\gamma.$$
On the other hand, by (\ref{beta}), we can select $z_n\in B(x_0,\beta);$ $v_n\in F(z_n)$ such that
$$ \|\lambda\bar y+\gamma u_n-(v_n-y_0)\| <\delta_n.$$
Therefore,
$$ \langle y^{1*}_n+y^{2*}_n,v_n-y_0\rangle<\langle y^{1*}_n+y^{2*}_n,\lambda\bar y+\gamma u_n\rangle+\delta_n<2\delta_n -(1-\delta_n)\gamma.$$
Since $(x^*_n,-y^{1*}_n,-y^{2*}_n)\in N(\gph G,(x_n,y^1_n,y_n^2)),$ then
$$\liminf_{n\to\infty} (\langle x^*_n, z_n-x_n\rangle -\langle y^{1*}_n+y^{2*}_n, v_n-y_0\rangle)\le 0. $$
Consequently,
$$ \liminf_{n\to\infty}\|x^*_n\|  \|z_n-x_n\|\ge\liminf_{n\to\infty}\langle y^{1*}_n+y^{2*}_n, v_n-y_0\rangle\ge (1-\delta_n)\gamma-2\delta_n.$$
By letting $n\to\infty,$ one obtains $\liminf_{n\to\infty}\|x^*_n\|\ge \gamma/\alpha,$ which shows (\ref{CC}).
\vskip 0.2cm
Suppose now that there exist $\tau>0,$ $\delta>0$ such that
$$ d(x,F^{-1}(y)) \le \tau d(y, F(x))\;\;\mbox{for all}\;\; (x,y)\in B(x_0,\delta)\times B(y_0,\delta),\; y\in F(x)+\cone B(\bar y,\delta).$$
In particular, one has
$$ F^{-1}(y) \not=\emptyset\quad\mbox{for all}\;\; y\in B(y_0,\varepsilon)\cap (y_0+\cone B(\bar y,\delta)).$$
Let $\varepsilon>0$ $\alpha>0$ be sufficiently small such that $\varepsilon<\delta\lambda$ and $\|\lambda\bar y\|+\varepsilon<\delta.$ Then,  for all $u\in B_Y$
$$z:=y_0+\lambda\bar y+\varepsilon u\in B(y_0,\delta)\cap (y_0+\cone B(\bar y,\delta)).$$
Hence, $F^{-1}(z)\not=\emptyset.$  It follows that $B(\lambda\bar y,\varepsilon)\subseteq F(X)-y_0,$ and the proof is complete.
\hfill\Halmos 

\section{Robustness of directional metric regularity}
The  characterizations of  directional metric regularity established in Theorem \ref{Direct-Char-strong-slope}, enable us  to derive  the following result on the stability of  directional metric regularity under perturbation. This result  has been first  obtained in \cite{ArAvIz07} under the inner semicontinuity  assumption. Then, when the image space $Y$ is a Banach space, Ioffe in \cite{RefIoffe} has extended  this    stability result (without   the inner semicontinuityassumption) with estimates sharper than the one in \cite{ArAvIz07}. Here, based on the mentioned characterizations, we prove this result for  which, the completeness of $Y$ is not necessary. 
\vskip 0.2cm
\begin{theorem} \label{Robust} Let $X$ be a complete metric space and $Y$ be a normed space. Let $F:X \rightrightarrows Y$ be a closed multifunction and $(x_0,y_0)\in \gph F.$ Suppose that $F$ is  metrically  regular with a modulus $\tau>0$ in the  direction $\bar y\in Y, $ i.e., there  exist  $\varepsilon>0,$ $\delta>0$ such that
\begin{equation}\label{DM1}
d(x,F^{-1}(y)) \le \tau d(y,F(x))\;\;\mbox{for all}\;(x,y)\in B((x_0,y_0),\varepsilon)\cap V(\bar y,\delta)\;\;\mbox {with}\; d(y,F(x))<\varepsilon.
\end{equation}
Let a  mapping $g: X\rightarrow Y$ be locally Lipschitz around $x_0$ with a Lipschitz constant $L>0.$ Then $F+g$ is  metrically  regular in the  direction $\bar y$ at $(x_0,y_0+g(x_0))$ with 
$${\rm reg}_{\bar y}(F+g) (x_0,y_0+g(x_0))\le \left(\frac{1-\gamma}{\tau(1+\gamma)}-L\right)^{-1},$$
provided 
$$ \alpha\in (0,1),\;\;\gamma:=\frac{\alpha\|\bar y\|}{\|\bar y\|+\delta(1-\alpha)};\;\; L<\frac{\delta(1-\alpha)\alpha}{\tau((1+\alpha)\|\bar y\|+\delta(1-\alpha))}.$$
\end{theorem}
\vskip 0.2cm
\proof{Proof.}  Let $\varepsilon,\delta,\gamma, \alpha, L$ as in Theorem \ref{Robust}. Let $g:X\to Y$ be Lipschitz with constant $L$ on $B(x_0,\varepsilon).$ To simplify the notations, denote by
$$ V_F(\delta):= \{(x,y):\;\; y\in F(x)+\cone B(\bar y,\delta) \};\;\;V_{(F+g)}(\delta):= \{(x,y):\;\; y\in F(x)+g(x)+\cone B(\bar y,\delta) \};$$
$$ V_{F,y}(\delta):= \{x\in X:\;\; y\in F(x)+\cone B(\bar y,\delta) \};\;\;V_{F+g,y}(\delta):= \{x\in X:\;\; y\in F(x)+g(x)+\cone B(\bar y,\delta) \},$$
and $\varphi_{V_F}(x,y),$ (resp. $\varphi_{V_{F+g}}(x,y)$) the lower semicontinuous envelope relative to $V_F$ (resp. $V_{F+g}$) of $F$ (resp. $F+g$).
Obviously, 
$$ \varphi_{V_{F+g}}(x,y)= \varphi_{V_F}(x,y-g(x)),\;\;\mbox{for all}\;\; (x,y)\in X\times Y.$$
According to Theorem \ref{Direct-Char-strong-slope} (ii), it suffices to prove that 
\begin{equation}
|\Gamma\varphi_{V_{F+g}}(\cdot,y)|(x)\ge \left(\frac{1-\gamma}{\tau(1+\gamma)}-L\right),
\end{equation}
whenever 
\begin{equation}\label{TT}(x,y) \in B( (x_0,y_0+g(x_0),\eta)\;\;\mbox{satisfies}\;\;x\in \cl V_{F+g,y}(\rho);\;\; d(y, F(x)+g(x))<\eta,\end{equation}
where $\rho:=\delta(1-\alpha)$ and $\eta=\min\{\varepsilon/(L+2),\varepsilon/(8\tau)\}$.
\vskip 0.2cm
Let $x,y$ be as in (\ref{TT}). Then   select sequences  $(\lambda_n)_{n\in\N}, (z_n)_{n\in\N}, (x_n))_{n\in\N}$ satisfying  $\lambda_n>0, z_n\in B_X, (x_n)\to x$  and such that
\begin{equation}\label{golf}
y-g(x_n)\in F(x_n) +\lambda_n(\bar y+\rho z_n),\quad \lim_{n\to\infty}d(y,F(x_n) +g(x_n))=\varphi_{V_{F+g}(\rho)}(x,y).\end{equation}
Note that since $(x_n)$ tends to $x$ and $x\in B(x_0,\eta)$, then for $n$ large we have 
\begin{equation}\label{golf1}
d(y, F(x_n)+ g(x_n))<\eta.
\end{equation}
 Setting
\begin{equation}\label{pluie}
 t_n:=\alpha \varphi_{V_{F+g}(\rho)}(x_n,y)/(\|\bar y\|+\rho),\end{equation} we observe that 
 \begin{equation}\label{porcelaine}t_n\|\bar y\| <\varphi_{V_{F+g}(\rho)}(x_n,y)<\eta,\end{equation}
and
$$d(y, F(x_n)+g(x_n))\leq \lambda_n\Vert \bar y +\rho z_n\Vert\leq \lambda_n(\Vert \bar y \Vert +\rho)\;\text{for some }\; z_n \;\text{with}\; \vert z_n\Vert = 1.$$
This yields, 
$$t_n(\|\bar y\|+\rho)/\alpha\leq\varphi_{V_{F+g}(\rho)}(x_n,y)\leq  d(y, F(x_n)+g(x_n))\leq\lambda_n(\Vert \bar y \Vert +\rho). $$
Consequently,   
\begin{equation}\label{golf2}
t_n/ \lambda_n\leq \alpha. \end{equation}
Observe also that 
$$
\begin{array}{ll}
 \lambda_n(\bar y+\rho z_n)-t_n\bar y\\
 &=  (\lambda_n-t_n)\bar y+\lambda_n\rho z_n\\
& = (\lambda_n-t_n)(\bar y +\frac{\lambda_n\rho}{\lambda_n-t_n} z_n).\\
\end{array}
$$
According   to (\ref{golf2}) 
$$ \frac{\lambda_n-t_n}{\lambda_n} = 1-\frac{t_n}{\lambda_n}\geq 1-\alpha$$
and therefore
$$\frac{\lambda_n}{\lambda_n-t_n}\rho\leq \frac{\rho}{1-\alpha} = \delta.$$
Hence, $\lambda_n(\bar y+\rho z_n)-t_n\bar y\in \cone B(\bar y,\varepsilon)$ and 
thanks to  (\ref{golf}), this yields
\begin{equation}\label{golf4}
y-g(x_n) -t_n \bar y \in F(x_n)+\cone B(\bar y , \delta).
\end{equation}
Moreover,
\begin{equation}\label{10} \|y-g(x_n)-t_n\bar y-y_0\| \le \|y-g(x_0)-y_0\| +\|g(x_n)-g(x_0)\| +t_n\|\bar y\|<(2+L)\eta=\varepsilon;\end{equation}
and combining \ref{golf1})  and (\ref{porcelaine})  we also have
\begin{equation}\label{11} d(y-g(x_n)-t_n\bar y, F(x_n)) \le d(y-g(x_n), F(x_n)) +t_n\|\bar y\|<2\eta<\frac{2\varepsilon}{L+2}<\varepsilon.\end{equation}
From (\ref{10}) and (\ref{11}) we deduce that  
\begin{itemize}
\item $y-g(x_n)-t_n\bar y\in B(y_0,\varepsilon);$
\item $d(y-g(x_n)-t_n\bar y, F(x_n))<\varepsilon;$
\item $(x_n,y-g(x_n)-t_n \bar y)\in V_F(\delta).$
\end{itemize}
Hence according to  
Proposition 4 (ii)
 we have 
\begin{align}
d(x_n, F^{-1} &(y-g(x_n)-t_n \bar y))\nonumber \\&<\tau\varphi_{V_{F(\rho)}}(x_n, y-g(x_n)-t_n\bar y)\nonumber \\
&\le \tau(\varphi_{V_{F+g}(\rho)}(x_n,y)+t_n\|\bar y\|)\nonumber\\
&=\tau t_n\frac{\Big( (1+\alpha)\|\bar y\|+\rho\Big)}{\alpha}\qquad \text{thanks to (\ref{pluie})}.\label{Ex-bis}\end{align}
Using the fact that  $t_n\Vert \bar y\Vert<\eta $ and $\varphi_{V_{F+g}(\rho)}(x_n,y)\le d(y-g(x_n), F(x_n)<\eta$, we obtain
$$d(x_n, F^{-1} (y-g(x_n)-t_n \bar y)<2\tau \eta.$$
By the choice of $\eta$, we derive $d(x_n, F^{-1} (y-g(x_n)-t_n \bar y)<\varepsilon /2$, and therefore for any $r\in (0,1),$  the existence of some $u_n\in F^{-1} (y-g(x_n)-t_n)$ such that 
$$d(x_n, u_n)<\tau (1+r) t_n( (1+\alpha)\|\bar y\|+\rho)/\alpha<\varepsilon /2.$$ 
Since $(x_n)\to x\in B(x_0,\eta)$, for $n$ sufficiently large we have 
$d(x_n, x_0)\leq d(x_n,x)+d(x,x_0)<\varepsilon /2 +\eta <\varepsilon,$  so that $u_n\in B(x_0, \varepsilon).$
Since $u_n\in F^{-1}(y-g(x_n)-t_n\bar y)\cap B(x_0, \varepsilon)$ and  by the Lipschitz property of  $g$ on $B(x_0,\varepsilon):$
$$\Vert  g(u_n)-g(x_n)\Vert \leq L d(u_n,x_n),$$ then 
$$ y \in F(u_n)+g(x_n)+t_n\bar y\subseteq F(u_n)+g(u_n)+t_n\Big(\bar y+L \frac{d(u_n,x_n)}{t_n}B_Y\Big). $$
By the definition of $L,$ for $r$ sufficiently  small, one obtains
$$ y\in F(u_n) +g(u_n)+\cone B(\bar y,\rho). $$
Therefore,
\begin{equation} \label{Ex-bis2}
\varphi_{V_{F+g}(\rho)}(u_n,y) \le d(y-g(u_n),F(u_n))\le t_n\|\bar y\|+Ld(x_n,u_n).
\end{equation}
As $t_n\|\bar y\| \le\alpha\varphi_{V_{F+g}(\rho)}(x_n,y)$ with $\alpha\in (0,1),$
it  follows that $\liminf_{n\to\infty}d(x_n,u_n)>0.$  Therefore, one has
\begin{align*} 
&\displaystyle\liminf_{n\to\infty}\frac{\varphi_{V_{F+g}(\rho)}(x,y)\varphi_{V_{F+g}(\rho)}(u_n,y)}{d(x,u_n)}
\\&=\displaystyle\liminf_{n\to\infty}\frac{\varphi_{V_{F+g}(\rho)}(x_n,y)-\varphi_{V_{F+g}(\rho)}(u_n,y)}{d(x_n,u_n)}\\&\ge\displaystyle\liminf_{n\to\infty} \frac{t_n(\|\bar y\|+\rho)/\alpha-t_n\|\bar y\|}{t_n((1+r)(\|\bar y\|+\rho)/\alpha+t_n\|\bar y\|}-L\\&=\displaystyle\frac{(\|\bar y\|+\rho)/\alpha-\|\bar y\|}{((1+r)(\|\bar y\|+\rho)/\alpha+\|\bar y\|}-L.
\end{align*}
As $r>0$ is  arbitrary small, one obtains
$$ |\Gamma\varphi_{V_{F+g}(\rho)}(\cdot,y)|(x) \ge\frac{1-\gamma}{\tau(1+\gamma)}-L,$$
which completes the proof.\hfill\Halmos 



%
%
%

\section*{Acknowledgments.}

Research  of Huynh Van Ngai was supported   by  VIASM (Vietnam Institute of Avanced Study on Mathematics).

Research of Michel Th\'era was partially supported by  by  the Australian Research Council under grant DP-$110102011$ and by LIA ``FormathVietnam". 



\bibliographystyle{ormsv080}
\bibliography{myrefs.THERA}

\def\Dbar{\leavevmode\lower.6ex\hbox to 0pt{\hskip-.23ex \accent"16\hss}D}
  \def\cfac#1{\ifmmode\setbox7\hbox{$\accent"5E#1$}\else
  \setbox7\hbox{\accent"5E#1}\penalty 10000\relax\fi\raise 1\ht7
  \hbox{\lower1.15ex\hbox to 1\wd7{\hss\accent"13\hss}}\penalty 10000
  \hskip-1\wd7\penalty 10000\box7}
\begin{thebibliography}{60}
\expandafter\ifx\csname natexlab\endcsname\relax\def\natexlab#1{#1}\fi
\expandafter\ifx\csname url\endcsname\relax
  \def\url#1{{\tt #1}}\fi
\expandafter\ifx\csname urlprefix\endcsname\relax\def\urlprefix{URL }\fi
\expandafter\ifx\csname urlstyle\endcsname\relax
  \expandafter\ifx\csname doi\endcsname\relax
  \def\doi#1{doi:\discretionary{}{}{}#1}\fi \else
  \expandafter\ifx\csname doi\endcsname\relax
  \def\doi{doi:\discretionary{}{}{}\begingroup \urlstyle{rm}\Url}\fi \fi

\bibitem[{Arutyunov et~al.(2007)Arutyunov, Avakov, and Izmailov}]{ArAvIz07}
Arutyunov, A.V., E.R. Avakov, A.F. Izmailov. 2007.
\newblock Directional regularity and metric regularity.
\newblock {\it SIAM J. Optim.\/} {\bf 18}(3) 810--833.
\newblock \doi{10.1137/060651616}.

\bibitem[{Az{\'e}.(2003)}]{Refaze}
Az{\'e}., D. 2003.
\newblock A survey on error bounds for lower semicontinuous functions.
\newblock {\it Proceedings of 2003 {MODE}-{SMAI} {C}onference\/}, {\it ESAIM
  Proc.\/}, vol.~13. EDP Sci., Les Ulis, 1--17.

\bibitem[{Az{\'e}(2006)}]{Aze06}
Az{\'e}, D. 2006.
\newblock A unified theory for metric regularity of multifunctions.
\newblock {\it J. Convex Anal.\/} {\bf 13}(2) 225--252.

\bibitem[{Az{\'e} and Corvellec(2004)}]{RefAC2}
Az{\'e}, D., J.-N. Corvellec. 2004.
\newblock Characterizations of error bounds for lower semicontinuous functions
  on metric spaces.
\newblock {\it ESAIM Control Optim. Calc. Var.\/} {\bf 10}(3) 409--425.
\newblock \doi{10.1051/cocv:2004013}.

\bibitem[{Bonnans and Shapiro(2000)}]{RefBonS}
Bonnans, J.~F., A.~Shapiro. 2000.
\newblock {\it Perturbation analysis of optimization problems\/}.
\newblock Springer Series in Operations Research, Springer-Verlag, New York.

\bibitem[{Borwein and Zhu(2005)}]{BorweinZhuBook}
Borwein, J.~M., Q.J. Zhu. 2005.
\newblock {\it Techniques of Variational Analysis\/}.
\newblock CMS Books in Mathematics/Ouvrages de Math\'ematiques de la SMC, 20,
  Springer-Verlag, New York.

\bibitem[{Borwein and Zhuang(1988)}]{BorZhuang88}
Borwein, J.~M., D.~M. Zhuang. 1988.
\newblock Verifiable necessary and sufficient conditions for openness and
  regularity for set-valued and single-valued mapps.
\newblock {\it J. Math. Anal. Appl.\/} {\bf 134} 441--459.

\bibitem[{Borwein and Dontchev(2003)}]{BD}
Borwein, J.M., A.L. Dontchev. 2003.
\newblock On the {B}artle-{G}raves theorem.
\newblock {\it Proc. Amer. Math. Soc.\/} {\bf 131}(8) 2553--2560.
\newblock \doi{10.1090/S0002-9939-03-07229-0}.

\bibitem[{Borwein and Zhu(1996)}]{BorZhu96}
Borwein, J.M., Q.J. Zhu. 1996.
\newblock Viscosity solutions and viscosity subderivatives in smooth {B}anach
  spaces with applications to metric regularity.
\newblock {\it SIAM J. Contr. Optim.\/} {\bf 34} 1568--1591.

\bibitem[{Bosch et~al.(2004)Bosch, Jourani, and Henrion}]{BJH}
Bosch, P., A.~Jourani, R.~Henrion. 2004.
\newblock Sufficient conditions for error bounds and applications.
\newblock {\it Appl. Math. Optim.\/} {\bf 50}(2) 161--181.

\bibitem[{Clarke(1983)}]{Clarke83}
Clarke, F.H. 1983.
\newblock {\it Optimization and Nonsmooth Analysis\/}.
\newblock Canadian Mathematical Society Series of Monographs and Advanced
  Texts, John Wiley \& Sons Inc., New York.
\newblock A Wiley-Interscience Publication.

\bibitem[{Cominetti(1990)}]{RefCom}
Cominetti, R. 1990.
\newblock Metric regularity, tangent sets, and second-order optimality
  conditions.
\newblock {\it Appl. Math. Optim.\/} {\bf 21}(3) 265--287.
\newblock \doi{10.1007/BF01445166}.

\bibitem[{Conway(1990)}]{conway}
Conway, J.B. 1990.
\newblock {\it A course in functional analysis\/}, {\it Graduate Texts in
  Mathematics\/}, vol.~96.
\newblock 2nd ed. Springer-Verlag, New York.

\bibitem[{DeGiorgi et~al.(1980)DeGiorgi, Marino, and Tosques}]{RefDMT}
DeGiorgi, E., A.~Marino, M.~Tosques. 1980.
\newblock Problems of evolution in metric spaces and maximal decreasing curve.
\newblock {\it Atti Accad. Naz. Lincei Rend. Cl. Sci. Fis. Mat. Natur. (8)\/}
  {\bf 68}(3) 180--187.

\bibitem[{Dmitruk et~al.(1980)Dmitruk, Milyutin, and Osmolovsky}]{RefDMO80}
Dmitruk, A.~V., A.A. Milyutin, N.P. Osmolovsky. 1980.
\newblock Lyusternik's theorem and the theory of extrema.
\newblock {\it Uspekhy Mat. Nauk\/} {\bf 35} 11--46.
\newblock In Russian.

\bibitem[{Dmitruk and Kruger(2008)}]{DmiKru08.1}
Dmitruk, A.V., A.~Y. Kruger. 2008.
\newblock Metric regularity and systems of generalized equations.
\newblock {\it J. Math. Anal. Appl.\/} {\bf 342}(2) 864--873.
\newblock \doi{10.1016/j.jmaa.2007.12.057}.

\bibitem[{Dmitruk and Kruger(2009)}]{DmiKru09.1}
Dmitruk, A.V., A.~Y. Kruger. 2009.
\newblock Extensions of metric regularity.
\newblock {\it Optimization\/} {\bf 58}(5) 561--584.
\newblock \doi{10.1080/02331930902928674}.

\bibitem[{Dontchev(1996)}]{dontchev}
Dontchev, A.L. 1996.
\newblock The {G}raves theorem revisited.
\newblock {\it J. Convex Anal.\/} {\bf 3}(1) 45--53.

\bibitem[{Dontchev. et~al.(2003)Dontchev., Lewis, and
  Rockafellar}]{DonLewRock03}
Dontchev., A.L., A.S. Lewis, R.T. Rockafellar. 2003.
\newblock The radius of metric regularity.
\newblock {\it Trans. Amer. Math. Soc.\/} {\bf 355}(2) 493--517.
\newblock \doi{10.1090/S0002-9947-02-03088-X}.

\bibitem[{Dontchev and Rockafellar(2009)}]{DonRoc09}
Dontchev, A.L., R.T. Rockafellar. 2009.
\newblock {\it Implicit Functions and Solution Mappings. A View from
  Variational Analysis\/}.
\newblock Springer Monographs in Mathematics, Springer, Dordrecht.

\bibitem[{Ekeland(1974)}]{RefE}
Ekeland, I. 1974.
\newblock On the variational principle.
\newblock {\it J. Math. Anal. Appl.\/} {\bf 47} 324--353.

\bibitem[{Fabian et~al.(2010)Fabian, Henrion, Kruger, and
  Outrata}]{FabHenKruOut10}
Fabian, M.~J., R.~Henrion, A.Y. Kruger, J.V. Outrata. 2010.
\newblock Error bounds: necessary and sufficient conditions.
\newblock {\it Set-Valued Var. Anal.\/} {\bf 18}(2) 121--149.
\newblock \doi{10.1007/s11228-010-0133-0}.

\bibitem[{Fabian et~al.(2012)Fabian, Henrion, Kruger, and
  Outrata}]{FabHenKruOut12}
Fabian, M.J., R.~Henrion, A.Y. Kruger, J.~V. Outrata. 2012.
\newblock About error bounds in metric spaces.
\newblock Diethard Klatte, Hans-Jacob L\"uthi, Karl Schmedders, eds., {\it
  Operations Research Proceedings 2011. Selected papers of the Int. Conf.
  Operations Research (OR 2011), August 30 -- September 2, 2011, Zurich,
  Switzerland\/}. Springer-Verlag, Berlin, 33--38.

\bibitem[{Giannessi(2005)}]{Gia05}
Giannessi, F. 2005.
\newblock {\it Constrained Optimization and Image Space Analysis. {V}ol. 1:
  Separation of Sets and Optimality Conditions\/}, {\it Mathematical Concepts
  and Methods in Science and Engineering\/}, vol.~49.
\newblock Springer, New York.

\bibitem[{Graves(1950)}]{RefGr}
Graves, L.M. 1950.
\newblock Some mapping theorems.
\newblock {\it Duke Math. J.\/} {\bf 17} 111--114.

\bibitem[{Hoffman(1952)}]{Hof52}
Hoffman, A.J. 1952.
\newblock On approximate solutions of systems of linear inequalities.
\newblock {\it J. Research Nat. Bur. Standards\/} {\bf 49} 263--265.

\bibitem[{Huynh et~al.(2010)Huynh, Kruger, and Th{\'e}ra}]{NgaiKruThe10}
Huynh, V.N., A.~Y. Kruger, M.~Th{\'e}ra. 2010.
\newblock Stability of error bounds for semi-infinite convex constraint
  systems.
\newblock {\it SIAM J. Optim.\/} {\bf 20}(4) 2080--2096.
\newblock \doi{10.1137/090767819}.

\bibitem[{Huynh et~al.(2012)Huynh, Kruger, and Th{\'e}ra}]{NgaKruThe12}
Huynh, V.N., A.Y. Kruger, M.~Th{\'e}ra. 2012.
\newblock Slopes of multifunctions and extensions of metric regularity.
\newblock {\it Vietnam J. Math.\/} {\bf 40}(2-3) 355--369.

\bibitem[{Huynh et~al.(2013)Huynh, Nguyen, and Th{\'e}ra}]{NTT}
Huynh, V.N., H.T. Nguyen, M.~Th{\'e}ra. 2013.
\newblock Implicit multifunctions theorems in complete metric spaces.
\newblock {\it Math. Program.\/} \doi{10.1007/s10107-013-0673-9}.

\bibitem[{Huynh and Th{\'e}ra(2001)}]{RefNT1}
Huynh, V.N., M.~Th{\'e}ra. 2001.
\newblock Metric inequality, subdifferential calculus and applications.
\newblock {\it Set-Valued Anal.\/} {\bf 9}(1-2) 187--216.
\newblock \doi{10.1023/A:1011291608129}.
\newblock Wellposedness in optimization and related topics (Gargnano, 1999).

\bibitem[{Huynh and Th{\'e}ra(2004)}]{RefNT3}
Huynh, V.N., M.~Th{\'e}ra. 2004.
\newblock Error bounds and implicit multifunction theorem in smooth {B}anach
  spaces and applications to optimization.
\newblock {\it Set-Valued Anal.\/} {\bf 12}(1-2) 195--223.
\newblock \doi{10.1023/B:SVAN.0000023396.58424.98}.

\bibitem[{Huynh and Th{\'e}ra(2005)}]{NgaiThe05}
Huynh, V.N., M.~Th{\'e}ra. 2005.
\newblock Error bounds for convex differentiable inequality systems in {B}anach
  spaces.
\newblock {\it Math. Program., Ser. B\/} {\bf 104}(2-3) 465--482.

\bibitem[{Huynh and Th{\'e}ra(2008)}]{NgaiThe08}
Huynh, V.N., M.~Th{\'e}ra. 2008.
\newblock Error bounds in metric spaces and application to the perturbation
  stability of metric regularity.
\newblock {\it SIAM J. Optim.\/} {\bf 19}(1) 1--20.

\bibitem[{Huynh and Th{\'e}ra(2009)}]{NgaiThe09}
Huynh, V.N., M.~Th{\'e}ra. 2009.
\newblock Error bounds for systems of lower semicontinuous functions in
  {A}splund spaces.
\newblock {\it Math. Program., Ser. B\/} {\bf 116}(1-2) 397--427.

\bibitem[{Ioffe(1979)}]{RefIo3}
Ioffe, A.D. 1979.
\newblock Regular points of {L}ipschitz functions.
\newblock {\it Trans. Amer. Math. Soc.\/} {\bf 251} 61--69.
\newblock \doi{10.2307/1998683}.

\bibitem[{Ioffe(2000)}]{Io00}
Ioffe, A.D. 2000.
\newblock Metric regularity and subdifferential calculus.
\newblock {\it Uspekhi Mat. Nauk\/} {\bf 55}(3(333)) 103--162.
\newblock In Russian.

\bibitem[{Ioffe(2001)}]{RefIo2}
Ioffe, A.D. 2001.
\newblock Towards metric theory of metric regularity.
\newblock {\it Approximation, optimization and mathematical economics
  ({P}ointe-\`a-{P}itre, 1999)\/}. Physica, Heidelberg, 165--176.

\bibitem[{Ioffe(2010)}]{RefIoffe}
Ioffe, A.D. 2010.
\newblock On regularity concepts in variational analysis.
\newblock {\it J. Fixed Point Theory Appl.\/} {\bf 8}(2) 339--363.
\newblock \doi{10.1007/s11784-010-0021-0}.

\bibitem[{Ioffe(2013)}]{ioffe2013}
Ioffe, A.D. 2013.
\newblock Convexity and variational analysis.
\newblock {\it Math. Program.\/} {\bf To appear}.

\bibitem[{Ioffe and Outrata(2008)}]{IofOut08}
Ioffe, A.D., J.V. Outrata. 2008.
\newblock On metric and calmness qualification conditions in subdifferential
  calculus.
\newblock {\it Set-Valued Anal.\/} {\bf 16}(2-3) 199--227.

\bibitem[{Jourani(2000)}]{RefJ}
Jourani, A. 2000.
\newblock Hoffman's error bound, local controllability, and sensitivity
  analysis.
\newblock {\it SIAM J. Control Optim.\/} {\bf 38}(3) 947--970.
\newblock \doi{10.1137/S0363012998339216}.

\bibitem[{Jourani and Thibault(1995)}]{RefJT}
Jourani, A., L.~Thibault. 1995.
\newblock Metric regularity and subdifferential calculus in {B}anach spaces.
\newblock {\it Set-Valued Anal.\/} {\bf 3}(1) 87--100.
\newblock \doi{10.1007/BF01033643}.

\bibitem[{Jourani and Thibault(1999)}]{JT1}
Jourani, A., L.~Thibault. 1999.
\newblock Coderivatives of multivalued mappings, locally compact cones and
  metric regularity.
\newblock {\it Nonlinear Anal.\/} {\bf 35}(7, Ser. A: Theory Methods) 925--945.
\newblock \doi{10.1016/S0362-546X(98)00031-5}.

\bibitem[{Klatte and Kummer(2002{\natexlab{a}})}]{kummer2002}
Klatte, D., B.~Kummer. 2002{\natexlab{a}}.
\newblock {\it Nonsmooth equations in optimization\/}, {\it Nonconvex
  Optimization and its Applications\/}, vol.~60.
\newblock Kluwer Academic Publishers, Dordrecht.

\bibitem[{Klatte and Kummer(2002{\natexlab{b}})}]{KlaKum02}
Klatte, D., B.~Kummer. 2002{\natexlab{b}}.
\newblock {\it Nonsmooth Equations in Optimization. Regularity, Calculus,
  Methods and Applications\/}, {\it Nonconvex Optimization and its
  Applications\/}, vol.~60.
\newblock Kluwer Academic Publishers, Dordrecht.

\bibitem[{Kruger et~al.(2010)Kruger, Huynh, and Th{\'e}ra}]{KruNgaiThe10}
Kruger, A.~Y., V.N. Huynh, M.~Th{\'e}ra. 2010.
\newblock Stability of error bounds for convex constraint systems in {B}anach
  spaces.
\newblock {\it SIAM J. Optim.\/} {\bf 20}(6) 3280--3296.
\newblock \doi{10.1137/100782206}.

\bibitem[{Kummer(1999)}]{Kum99}
Kummer, B. 1999.
\newblock Metric regularity: {C}haracterizations, nonsmooth variations and
  successive approximation.
\newblock {\it Optimization\/} {\bf 46} 247--281.

\bibitem[{Lyusternik(1934)}]{Lyusternik}
Lyusternik, L.A. 1934.
\newblock On conditional extrema of functionals.
\newblock {\it Math. Sbornik\/} {\bf 41} 390--401.
\newblock In Russian.

\bibitem[{Mordukhovich(2006{\natexlab{a}})}]{Borisbook1}
Mordukhovich, B.~S. 2006{\natexlab{a}}.
\newblock {\it Variational analysis and generalized differentiation. {I}\/},
  {\it Grundlehren der Mathematischen Wissenschaften [Fundamental Principles of
  Mathematical Sciences]\/}, vol. 330.
\newblock Springer-Verlag, Berlin.
\newblock Basic theory.

\bibitem[{Mordukhovich(2006{\natexlab{b}})}]{Borisbook2}
Mordukhovich, B.~S. 2006{\natexlab{b}}.
\newblock {\it Variational analysis and generalized differentiation. {II}\/},
  {\it Grundlehren der Mathematischen Wissenschaften [Fundamental Principles of
  Mathematical Sciences]\/}, vol. 331.
\newblock Springer-Verlag, Berlin.
\newblock Applications.

\bibitem[{Mordukhovich and Outrata(2007)}]{Mo-Ou2007}
Mordukhovich, B.~S., J.V. Outrata. 2007.
\newblock Coderivative analysis of quasi-variational inequalities with
  applications to stability and optimization.
\newblock {\it SIAM J. Optim.\/} {\bf 18}(2) 389--412.
\newblock \doi{10.1137/060665609}.

\bibitem[{Mordukhovich and Shao(1997)}]{RefMorS}
Mordukhovich, B.~S., Y.~Shao. 1997.
\newblock Stability of set-valued mappings in infinite dimensions: point
  criteria and applications.
\newblock {\it SIAM J. Control Optim.\/} {\bf 35}(1) 285--314.
\newblock \doi{10.1137/S0363012994278171}.

\bibitem[{Mordukhovich and Wang(2004)}]{RefMorW}
Mordukhovich, B.~S., B.~Wang. 2004.
\newblock Restrictive metric regularity and generalized differential calculus
  in {B}anach spaces.
\newblock {\it Int. J. Math. Math. Sci.\/} (49-52) 2653--2680.
\newblock \doi{10.1155/S0161171204405183}.

\bibitem[{Outrata et~al.(1998)Outrata, Ko{\v{c}}vara, and Zowe}]{Outrata98}
Outrata, J., M.~Ko{\v{c}}vara, J.~Zowe. 1998.
\newblock {\it Nonsmooth approach to optimization problems with equilibrium
  constraints\/}, {\it Nonconvex Optimization and its Applications\/}, vol.~28.
\newblock Kluwer Academic Publishers, Dordrecht.
\newblock Theory, applications and numerical results.

\bibitem[{Penot(1989)}]{Pen89}
Penot, J.-P. 1989.
\newblock Metric regularity, openness and {L}ipschitz behavior of
  multifunctions.
\newblock {\it Nonlinear Anal.\/} {\bf 13} 629--643.

\bibitem[{Penot(2013)}]{penot2013}
Penot, J.-P. 2013.
\newblock {\it Calculus Without Derivatives\/}.
\newblock Graduate Texts in Mathematics, Springer Verlag, Berlin.

\bibitem[{Robinson(1973)}]{robinson73}
Robinson, S.M. 1973.
\newblock An inverse-function theorem for a class of multivalued functions.
\newblock {\it Proc. Amer. Math. Soc.\/} {\bf 41} 211--218.

\bibitem[{Robinson(1980)}]{robinson80}
Robinson, S.M. 1980.
\newblock Strongly regular generalized equations.
\newblock {\it Math. Oper. Res.\/} {\bf 5}(1) 43--62.

\bibitem[{Rockafellar and Wets(1998)}]{RockBook2002}
Rockafellar, R.T., R.J.-B. Wets. 1998.
\newblock {\it Variational analysis\/}, {\it Grundlehren der Mathematischen
  Wissenschaften [Fundamental Principles of Mathematical Sciences]\/}, vol.
  317.
\newblock Springer-Verlag, Berlin.

\bibitem[{Ursescu(1975)}]{urcescu}
Ursescu, C. 1975.
\newblock Multifunctions with convex closed graph.
\newblock {\it Czechoslovak Math. J.\/} {\bf 25(100)}(3) 438--441.

\end{thebibliography}
\end{document}